\documentclass[letterpaper,12pt]{amsart}
\usepackage{latexsym}
\usepackage{amssymb}
\usepackage{amsmath}
\usepackage{ctable}
\usepackage{qtree}
\usepackage{varioref}
\usepackage{enumerate}
\usepackage{pmat}
\usepackage{pstricks-add}

\begin{document}

\title{A catalog of Cayley-Dickson-like products}
\author{John W. Bales}
\date{June 2011}
\address{Department of Mathematics\\Tuskegee University\\Tuskegee, AL 36088\\   (June 2011)}
\email{jbales@mytu.tuskegee.edu}
\keywords{twisted group algebra, Cayley-Dickson algebra, quaternions, octonions, Fano plane}

\begin{abstract}
A catalog of all 32 Cayley-Dickson-like doubling products on ordered pairs $(a,b)\cdot(c,d)$ for which $(1,0)$ is the left and right identity and for
which $x\cdot x^*=x^*\cdot x=\|x\|^2$ given the conjugate $(a,b)^*=(a^*,-b).$ Only eight of these are true Cayley-Dickson doubling products, since 24
of them do not satisfy the quaternion properties. Each of the eight Cayley-Dickson products has a distinctive representation in the Fano Plane.
\end{abstract}

\maketitle

\section{Product of finite sequences of real numbers}

In this catalog, Cayley-Dickson-like products are considered to be products of finite sequences of real numbers. This may
seem odd since Cayley-Dickson products are normally considered to be products defined on ordered pairs. 
In order to define a Cayley-Dickson-like product of any two finite sequences of real numbers, every finite sequence $x$ will
be identified with an infinite sequence of the form

\begin{equation*}
  x=x_0,x_1,x_2,\cdots,x_{n-1},0,0,0,\cdots
\end{equation*}

and the ordered pair of any two sequences $x$ and $y$ (whether finite or infinite) will be identified with the shuffling
of the two sequences

\begin{equation*}
  (x,y)=x_0,y_0,x_1,y_1,x_2,y_2\cdots
\end{equation*}

If $a$ is a real number, then $a$ will be identified with the sequence $a,0,0,0,\cdots.$

The basis vectors $i_0,i_1,i_2,\cdots$  for this space are defined recursively.

\begin{eqnarray*}
  i_0&=&1\\
  i_{2k}&=&(i_k,0)\text{\ for\ }k\ge0\\
  i_{2k+1}&=&(0,i_k)\text{\ for\ }k\ge0
\end{eqnarray*}

Thus it follows that 

\begin{eqnarray*}
  i_0&=&1,0,0,0,\cdots\\
  i_1&=&0,1,0,0,\cdots\\
  i_2&=&0,0,1,0,\cdots\\
  \vdots
\end{eqnarray*}

Furthermore, an involution $^*$ is defined on the basis vectors\\

$  i_k^*=\begin{cases}
           i_k & \text{\ if\ } k=0\\
           -i_k &  \text{\ if\ } k>0
         \end{cases}$

and is extended to the entire space by the linearity property.\\

It then follows that for all infinite sequences $x$ and $y$

\begin{equation}
  (x,y)^*=(x^*,-y)
\end{equation}

The norm $\|x\|$ of a finite sequence $x$ is its Euclidean norm and the inner product $\langle x,y \rangle$ of two
finite sequences $x$ and $y$ is their Euclidean inner product. Thus

\begin{equation*}
  \|(x,y)\|^2=\|x\|^2+\|y\|^2=\left(\|x\|^2+\|y\|^2,0\right)
\end{equation*}

\section{Cayley-Dickson-like Product}

A Cayley-Dickson-like product on finite sequences $(a,b)$ and $(c,d)$ is a product satisfying the following four properties:

\begin{eqnarray}
  (1,0)(b,c)&=&(b,c)\\
  (a,b)(1,0)&=&(a,b)\\
  (a,b)(a^*,-b)&=&\left(\|a\|^2+\|b\|^2,0\right)\\
  (c^*,-d)(c,d)&=&\left(\|c\|^2+\|d\|^2,0\right)
\end{eqnarray}

It is `Cayley-Dickson-like' rather than 'Cayley-Dickson' since we do not require the quaternion property.

The quaternion property is a property of the basis vectors. The statement that the quaternion property is satisfied means that if $p,q$ and $r$ are
positive integers and if $i_pi_q=i_r$ then $i_qi_p=-i_r$ and $i_qi_r=i_p.$

The product of $(a,b)(c,d)$ can be defined 32 different ways to achieve the four `Cayley-Dickson-like' properties, thus guaranteeing that
for all finite sequences $x$,

\begin{eqnarray}
  1\cdot x&= &x\cdot 1=x\\
  xx^*&=&x^*x=\|x\|^2
\end{eqnarray}

Only eight of the 32 Cayley-Dickson-like products satisfy the quaternion properties. For all eight of these, $i_1=i.$ For four of them $i_2=j$ and
$i_3=k.$ For the other four, $i_2=k$ and $i_3=j.$

\section{The 32 ways}

For each of the 32 Cayley-Dickson-like products of finite sequences, the recursive multiplication table for the basis vectors
is given. The table is divided into distinct $2\times2$ blocks. These blocks occur either on the corner, left edge,
top edge, diagonal or interior of the multiplication table.

An examination of the tables reveals that for all 32 variations it is true that $i_1i_2+i_2i_1=0$. Using this as a
basis step, it can be shown by induction that for all 32 variations

\begin{equation}\label{eq:neg}
  i_qi_p=-i_pi_q \text{\ for\ }0\ne p\ne q\ne 0
\end{equation}

\newpage

\subsection{$P_{0}  :(a,b)(c,d)=(ca-b^*d,da^*+bc)$}  
For all $r,s.$\\

$\begin{array}{lll} P_{0}&i_{2s}&i_{2s+1}\\\toprule i_{2r}&(i_si_r,0)&(0,i_si_r^*)  \\\midrule i_{2r+1}  
&(0,i_ri_s)  
&(-i_r^*i_s,0)\\\midrule \end{array}$

\begin{enumerate}
  \item For $r=s=0$ (The corner)

$\begin{array}{lll} P_{0}&i_{0}&i_{1}\\\toprule i_{0}&1&i_1\\\midrule i_{1}  
&i_1
&-1\\\midrule 
\end{array}$

  \item For $r>s=0$ (The left edge)

$\begin{array}{lll} P_{0}&i_{0}&i_{1}\\\toprule i_{2r}&i_{2r}&-i_{2r+1}\\\midrule i_{2r+1}  
&i_{2r+1}
&i_{2r}\\\midrule 
\end{array}$

  \item For $s>r=0$ (The top edge)

$\begin{array}{lll} P_{0}&i_{2s}&i_{2s+1}\\\toprule i_{0}&i_{2s}&i_{2s+1}\\\midrule i_{1}  
&i_{2s+1}
&-i_{2s}\\\midrule 
\end{array}$

  \item For $r=s>0$ (The diagonal)

$\begin{array}{lll} P_{0}&i_{2r}&i_{2r+1}\\\toprule i_{2r}&-1&i_1\\\midrule i_{2r+1}  
&-i_1
&-1\\\midrule 
\end{array}$
  
  \item For $0\ne r\ne s\ne0$ (The interior)
  
$\begin{array}{lll} P_{0}&i_{2s}&i_{2s+1}\\\toprule i_{2r}&(i_{s}i_{r},0)&-(0,i_{s}i_{r})\\\midrule
i_{2r+1} 
&(0,i_{r}i_{s})
&(i_{r}i_{s},0)\\\midrule 
\end{array}$

\end{enumerate}

\newpage

\subsection{$P_{1}  :(a,b)(c,d)=(ca-db^*,da^*+bc)$}  
For all $r,s.$\\

$\begin{array}{lll} P_{1}&i_{2s}&i_{2s+1}\\\toprule i_{2r}&(i_si_r,0)&(0,i_si_r^*)  \\\midrule i_{2r+1}  
&(0,i_ri_s)  
&(-i_si_r^*,0)\\\midrule \end{array}$
  
\begin{enumerate}
  \item For $r=s=0$ (The corner)

$\begin{array}{lll} P_{1}&i_{0}&i_{1}\\\toprule i_{0}&1&i_1\\\midrule i_{1}  
&i_1
&-1\\\midrule 
\end{array}$

  \item For $r>s=0$ (The left edge)

$\begin{array}{lll} P_{1}&i_{0}&i_{1}\\\toprule i_{2r}&i_{2r}&-i_{2r+1}\\\midrule i_{2r+1}  
&i_{2r+1}
&i_{2r}\\\midrule 
\end{array}$

  \item For $s>r=0$ (The top edge)

$\begin{array}{lll} P_{1}&i_{2s}&i_{2s+1}\\\toprule i_{0}&i_{2s}&i_{2s+1}\\\midrule i_{1}  
&i_{2s+1}
&-i_{2s}\\\midrule 
\end{array}$

  \item For $r=s>0$ (The diagonal)

$\begin{array}{lll} P_{1}&i_{2r}&i_{2r+1}\\\toprule i_{2r}&-1&i_1\\\midrule i_{2r+1}  
&-i_1
&-1\\\midrule 
\end{array}$
  
  \item For $0\ne r\ne s\ne0$ (The interior)
  
$\begin{array}{lll} P_{1}&i_{2s}&i_{2s+1}\\\toprule i_{2r}&(i_{s}i_{r},0)&-(0,i_{s}i_{r})\\\midrule
i_{2r+1} 
&(0,i_{r}i_{s})
&(i_{s}i_{r},0)\\\midrule 
\end{array}$

\end{enumerate}

\newpage

\subsection{$P_{2}  :(a,b)(c,d)=(ca-b^*d,a^*d+cb)$}  
For all $r,s.$\\

$\begin{array}{lll} P_{2}&i_{2s}&i_{2s+1}\\\toprule i_{2r}&(i_si_r,0)&(0,i_r^*i_s)  \\\midrule i_{2r+1}  
&(0,i_si_r)  
&(-i_r^*i_s,0)\\\midrule \end{array}$

\begin{enumerate}
  \item For $r=s=0$ (The corner)

$\begin{array}{lll} P_{2}&i_{0}&i_{1}\\\toprule i_{0}&1&i_1\\\midrule i_{1}  
&i_1
&-1\\\midrule 
\end{array}$

  \item For $r>s=0$ (The left edge)

$\begin{array}{lll} P_{2}&i_{0}&i_{1}\\\toprule i_{2r}&i_{2r}&-i_{2r+1}\\\midrule i_{2r+1}  
&i_{2r+1}
&i_{2r}\\\midrule 
\end{array}$

  \item For $s>r=0$ (The top edge)

$\begin{array}{lll} P_{2}&i_{2s}&i_{2s+1}\\\toprule i_{0}&i_{2s}&i_{2s+1}\\\midrule i_{1}  
&i_{2s+1}
&-i_{2s}\\\midrule 
\end{array}$

  \item For $r=s>0$ (The diagonal)

$\begin{array}{lll} P_{2}&i_{2r}&i_{2r+1}\\\toprule i_{2r}&-1&i_1\\\midrule i_{2r+1}  
&-i_1
&-1\\\midrule 
\end{array}$
  
  \item For $0\ne r\ne s\ne0$ (The interior)
  
$\begin{array}{lll} P_{2}&i_{2s}&i_{2s+1}\\\toprule i_{2r}&(i_{s}i_{r},0)&-(0,i_{r}i_{s})\\\midrule
i_{2r+1} 
&(0,i_{s}i_{r})
&(i_{r}i_{s},0)\\\midrule 
\end{array}$

\end{enumerate}

\newpage

\subsection{$P_{3}  :(a,b)(c,d)=(ca-db^*,a^*d+cb)$}  
For all $r,s.$\\

$\begin{array}{lll} P_{3}&i_{2s}&i_{2s+1}\\\toprule i_{2r}&(i_si_r,0)&(0,i_r^*i_s)  \\\midrule i_{2r+1}  
&(0,i_si_r)  
&(-i_si_r^*,0)\\\midrule \end{array}$

\begin{enumerate}
  \item For $r=s=0$ (The corner)

$\begin{array}{lll} P_{3}&i_{0}&i_{1}\\\toprule i_{0}&1&i_1\\\midrule i_{1}  
&i_1
&-1\\\midrule 
\end{array}$

  \item For $r>s=0$ (The left edge)

$\begin{array}{lll} P_{3}&i_{0}&i_{1}\\\toprule i_{2r}&i_{2r}&-i_{2r+1}\\\midrule i_{2r+1}  
&i_{2r+1}
&i_{2r}\\\midrule 
\end{array}$

  \item For $s>r=0$ (The top edge)

$\begin{array}{lll} P_{3}&i_{2s}&i_{2s+1}\\\toprule i_{0}&i_{2s}&i_{2s+1}\\\midrule i_{1}  
&i_{2s+1}
&-i_{2s}\\\midrule 
\end{array}$

  \item For $r=s>0$ (The diagonal)

$\begin{array}{lll} P_{3}&i_{2r}&i_{2r+1}\\\toprule i_{2r}&-1&i_1\\\midrule i_{2r+1}  
&-i_1
&-1\\\midrule 
\end{array}$
  
  \item For $0\ne r\ne s\ne0$ (The interior)
  
$\begin{array}{lll} P_{3}&i_{2s}&i_{2s+1}\\\toprule i_{2r}&(i_{s}i_{r},0)&-(0,i_{r}i_{s})\\\midrule
i_{2r+1} 
&(0,i_{s}i_{r})
&(i_{s}i_{r},0)\\\midrule 
\end{array}$

\end{enumerate}

\newpage

     \subsection{$P_{4}  :(a,b)(c,d)=(ac-b^*d,da^*+bc)$}  
For all $r,s.$\\

$\begin{array}{lll} P_{4}&i_{2s}&i_{2s+1}\\\toprule i_{2r}&(i_ri_s,0)&(0,i_si_r^*)  \\\midrule i_{2r+1}  
&(0,i_ri_s)  
&(-i_r^*i_s,0)\\\midrule \end{array}$

\begin{enumerate}
  \item For $r=s=0$ (The corner)

$\begin{array}{lll} P_{4}&i_{0}&i_{1}\\\toprule i_{0}&1&i_1\\\midrule i_{1}  
&i_1
&-1\\\midrule 
\end{array}$

  \item For $r>s=0$ (The left edge)

$\begin{array}{lll} P_{4}&i_{0}&i_{1}\\\toprule i_{2r}&i_{2r}&-i_{2r+1}\\\midrule i_{2r+1}  
&i_{2r+1}
&i_{2r}\\\midrule 
\end{array}$

  \item For $s>r=0$ (The top edge)

$\begin{array}{lll} P_{4}&i_{2s}&i_{2s+1}\\\toprule i_{0}&i_{2s}&i_{2s+1}\\\midrule i_{1}  
&i_{2s+1}
&-i_{2s}\\\midrule 
\end{array}$

  \item For $r=s>0$ (The diagonal)

$\begin{array}{lll} P_{4}&i_{2r}&i_{2r+1}\\\toprule i_{2r}&-1&i_1\\\midrule i_{2r+1}  
&-i_1
&-1\\\midrule 
\end{array}$
  
  \item For $0\ne r\ne s\ne0$ (The interior)
  
$\begin{array}{lll} P_{4}&i_{2s}&i_{2s+1}\\\toprule i_{2r}&(i_{r}i_{s},0)&-(0,i_{s}i_{r})\\\midrule
i_{2r+1} 
&-(0,i_{r}i_{s})
&(i_{r}i_{s},0)\\\midrule 
\end{array}$

\end{enumerate}

\newpage

\subsection{$P_{5}  :(a,b)(c,d)=(ac-db^*,da^*+bc)$}  
For all $r,s.$\\

$\begin{array}{lll} P_{5}&i_{2s}&i_{2s+1}\\\toprule i_{2r}&(i_ri_s,0)&(0,i_si_r^*)  \\\midrule i_{2r+1}  
&(0,i_ri_s)  
&(-i_si_r^*,0)\\\midrule \end{array}$

\begin{enumerate}
  \item For $r=s=0$ (The corner)

$\begin{array}{lll} P_{5}&i_{0}&i_{1}\\\toprule i_{0}&1&i_1\\\midrule i_{1}  
&i_1
&-1\\\midrule 
\end{array}$

  \item For $r>s=0$ (The left edge)

$\begin{array}{lll} P_{5}&i_{0}&i_{1}\\\toprule i_{2r}&i_{2r}&-i_{2r+1}\\\midrule i_{2r+1}  
&i_{2r+1}
&i_{2r}\\\midrule 
\end{array}$

  \item For $s>r=0$ (The top edge)

$\begin{array}{lll} P_{5}&i_{2s}&i_{2s+1}\\\toprule i_{0}&i_{2s}&i_{2s+1}\\\midrule i_{1}  
&i_{2s+1}
&-i_{2s}\\\midrule 
\end{array}$

  \item For $r=s>0$ (The diagonal)

$\begin{array}{lll} P_{5}&i_{2r}&i_{2r+1}\\\toprule i_{2r}&-1&-i_1\\\midrule i_{2r+1}  
&i_1
&-1\\\midrule 
\end{array}$
  
  \item For $0\ne r\ne s\ne0$ (The interior)
  
$\begin{array}{lll} P_{5}&i_{2s}&i_{2s+1}\\\toprule i_{2r}&(i_{r}i_{s},0)&-(0,i_{s}i_{r})\\\midrule
i_{2r+1} 
&(0,i_{r}i_{s})
&(i_{s}i_{r},0)\\\midrule 
\end{array}$

\end{enumerate}

\newpage\subsection{$P_{6}  :(a,b)(c,d)=(ac-b^*d,a^*d+cb)$}  
For all $r,s.$\\

$\begin{array}{lll} P_{6}&i_{2s}&i_{2s+1}\\\toprule i_{2r}&(i_ri_s,0)&(0,i_r^*i_s)  \\\midrule i_{2r+1}  
&(0,i_si_r)  
&(-i_r^*i_s,0)\\\midrule \end{array}$

\begin{enumerate}
  \item For $r=s=0$ (The corner)

$\begin{array}{lll} P_{6}&i_{0}&i_{1}\\\toprule i_{0}&1&i_1\\\midrule i_{1}  
&i_1
&-1\\\midrule 
\end{array}$

  \item For $r>s=0$ (The left edge)

$\begin{array}{lll} P_{6}&i_{0}&i_{1}\\\toprule i_{2r}&i_{2r}&-i_{2r+1}\\\midrule i_{2r+1}  
&i_{2r+1}
&i_{2r}\\\midrule 
\end{array}$

  \item For $s>r=0$ (The top edge)

$\begin{array}{lll} P_{6}&i_{2s}&i_{2s+1}\\\toprule i_{0}&i_{2s}&i_{2s+1}\\\midrule i_{1}  
&i_{2s+1}
&-i_{2s}\\\midrule 
\end{array}$

  \item For $r=s>0$ (The diagonal)

$\begin{array}{lll} P_{6}&i_{2r}&i_{2r+1}\\\toprule i_{2r}&-1&i_1\\\midrule i_{2r+1}  
&-i_1
&-1\\\midrule 
\end{array}$
  
  \item For $0\ne r\ne s\ne0$ (The interior)
  
$\begin{array}{lll} P_{6}&i_{2s}&i_{2s+1}\\\toprule i_{2r}&(i_{r}i_{s},0)&-(0,i_{r}i_{s})\\\midrule
i_{2r+1} 
&(0,i_{s}i_{r})
&(i_{r}i_{s},0)\\\midrule 
\end{array}$
\end{enumerate}

\newpage\subsection{$P_{7}  :(a,b)(c,d)=(ac-db^*,a^*d+cb)$}  
For all $r,s.$\\

$\begin{array}{lll} P_{7}&i_{2s}&i_{2s+1}\\\toprule i_{2r}&(i_ri_s,0)&(0,i_r^*i_s)  \\\midrule i_{2r+1}  
&(0,i_si_r) 
&(-i_si_r^*,0)\\\midrule \end{array}$

\begin{enumerate}
  \item For $r=s=0$ (The corner)

$\begin{array}{lll} P_{7}&i_{0}&i_{1}\\\toprule i_{0}&1&i_1\\\midrule i_{1}  
&i_1
&-1\\\midrule 
\end{array}$

  \item For $r>s=0$ (The left edge)

$\begin{array}{lll} P_{7}&i_{0}&i_{1}\\\toprule i_{2r}&i_{2r}&-i_{2r+1}\\\midrule i_{2r+1}  
&i_{2r+1}
&i_{2r}\\\midrule 
\end{array}$

  \item For $s>r=0$ (The top edge)

$\begin{array}{lll} P_{7}&i_{2s}&i_{2s+1}\\\toprule i_{0}&i_{2s}&i_{2s+1}\\\midrule i_{1}  
&i_{2s+1}
&-i_{2s}\\\midrule 
\end{array}$

  \item For $r=s>0$ (The diagonal)

$\begin{array}{lll} P_{7}&i_{2r}&i_{2r+1}\\\toprule i_{2r}&-1&i_1\\\midrule i_{2r+1}  
&-i_1
&-1\\\midrule 
\end{array}$
  
  \item For $0\ne r\ne s\ne0$ (The interior)
  
$\begin{array}{lll} P_{7}&i_{2s}&i_{2s+1}\\\toprule i_{2r}&(i_{r}i_{s},0)&-(0,i_{r}i_{s})\\\midrule
i_{2r+1} 
&(0,i_{s}i_{r})
&(i_{s}i_{r},0)\\\midrule 
\end{array}$

\end{enumerate}

     \newpage\subsection{$P_{8}   :(a,b)(c,d)=(ca-bd^*,da^*+bc)$}  
For all $r,s.$\\

$\begin{array}{lll} P_{8}&i_{2s}&i_{2s+1}\\\toprule i_{2r}&(i_si_r,0)&(0,i_si_r^*)  \\\midrule i_{2r+1}  
&(0,i_ri_s)  
&(-i_ri_s^*,0)\\\midrule \end{array}$
  
\begin{enumerate}
  \item For $r=s=0$ (The corner)

$\begin{array}{lll} P_{8}&i_{0}&i_{1}\\\toprule i_{0}&1&i_1\\\midrule i_{1}  
&i_1
&-1\\\midrule 
\end{array}$

  \item For $r>s=0$ (The left edge)

$\begin{array}{lll} P_{8}&i_{0}&i_{1}\\\toprule i_{2r}&i_{2r}&-i_{2r+1}\\\midrule i_{2r+1}  
&i_{2r+1}
&-i_{2r}\\\midrule 
\end{array}$

  \item For $s>r=0$ (The top edge)

$\begin{array}{lll} P_{8}&i_{2s}&i_{2s+1}\\\toprule i_{0}&i_{2s}&i_{2s+1}\\\midrule i_{1}  
&i_{2s+1}
&i_{2s}\\\midrule 
\end{array}$

  \item For $r=s>0$ (The diagonal)

$\begin{array}{lll} P_{8}&i_{2r}&i_{2r+1}\\\toprule i_{2r}&-1&i_1\\\midrule i_{2r+1}  
&-i_1
&-1\\\midrule 
\end{array}$
  
  \item For $0\ne r\ne s\ne0$ (The interior)
  
$\begin{array}{lll} P_{8}&i_{2s}&i_{2s+1}\\\toprule i_{2r}&(i_{s}i_{r},0)&-(0,i_{s}i_{r})\\\midrule
i_{2r+1} 
&(0,i_{r}i_{s})
&(i_{r}i_{s},0)\\\midrule 
\end{array}$

\end{enumerate}

\newpage\subsection{$P_{9}   :(a,b)(c,d)=(ca-d^*b,da^*+bc)$}
       
  For all $r,s.$\\
  
$\begin{array}{lll} P_{9}&i_{2s}&i_{2s+1}\\\toprule i_{2r}&(i_si_r,0)&(0,i_si_r^*)  \\\midrule i_{2r+1}  
&(0,i_ri_s)  
&(-i_s^*i_r,0)\\\midrule \end{array}$

\begin{enumerate}
  \item For $r=s=0$ (The corner)

$\begin{array}{lll} P_{9}&i_{0}&i_{1}\\\toprule i_{0}&1&i_1\\\midrule i_{1}  
&i_1
&-1\\\midrule 
\end{array}$

  \item For $r>s=0$ (The left edge)

$\begin{array}{lll} P_{9}&i_{0}&i_{1}\\\toprule i_{2r}&i_{2r}&-i_{2r+1}\\\midrule i_{2r+1}  
&i_{2r+1}
&-i_{2r}\\\midrule 
\end{array}$

  \item For $s>r=0$ (The top edge)

$\begin{array}{lll} P_{9}&i_{2s}&i_{2s+1}\\\toprule i_{0}&i_{2s}&i_{2s+1}\\\midrule i_{1}  
&i_{2s+1}
&i_{2s}\\\midrule 
\end{array}$

  \item For $r=s>0$ (The diagonal)

$\begin{array}{lll} P_{9}&i_{2r}&i_{2r+1}\\\toprule i_{2r}&-1&i_1\\\midrule i_{2r+1}  
&-i_1
&-1\\\midrule 
\end{array}$
  
  \item For $0\ne r\ne s\ne0$ (The interior)
  
$\begin{array}{lll} P_{9}&i_{2s}&i_{2s+1}\\\toprule i_{2r}&(i_{s}i_{s},0)&-(0,i_{s}i_{r})\\\midrule
i_{2r+1} 
&(0,i_{r}i_{s})
&(i_{s}i_{r},0)\\\midrule 
\end{array}$

\end{enumerate}

\newpage\subsection{$P_{10}   :(a,b)(c,d)=(ca-bd^*,a^*d+cb)$}  
For all $r,s.$\\

$\begin{array}{lll} P_{10}&i_{2s}&i_{2s+1}\\\toprule i_{2r}&(i_si_r,0)&(0,i_r^*i_s)  \\\midrule i_{2r+1}  
&(0,i_si_r)  
&(-i_ri_s^*,0)\\\midrule \end{array}$
  
\begin{enumerate}
  \item For $r=s=0$ (The corner)

$\begin{array}{lll} P_{10}&i_{0}&i_{1}\\\toprule i_{0}&1&i_1\\\midrule i_{1}  
&i_1
&-1\\\midrule 
\end{array}$

  \item For $r>s=0$ (The left edge)

$\begin{array}{lll} P_{10}&i_{0}&i_{1}\\\toprule i_{2r}&i_{2r}&-i_{2r+1}\\\midrule i_{2r+1}  
&i_{2r+1}
&-i_{2r}\\\midrule 
\end{array}$

  \item For $s>r=0$ (The top edge)

$\begin{array}{lll} P_{10}&i_{2s}&i_{2s+1}\\\toprule i_{0}&i_{2s}&i_{2s+1}\\\midrule i_{1}  
&i_{2s+1}
&i_{2s}\\\midrule 
\end{array}$

  \item For $r=s>0$ (The diagonal)

$\begin{array}{lll} P_{10}&i_{2r}&i_{2r+1}\\\toprule i_{2r}&-1&i_1\\\midrule i_{2r+1}  
&-i_1
&-1\\\midrule 
\end{array}$
  
  \item For $0\ne r\ne s\ne0$ (The interior)
  
$\begin{array}{lll} P_{10}&i_{2s}&i_{2s+1}\\\toprule i_{2r}&(i_{s}i_{r},0)&-(0,i_{r}i_{s})\\\midrule
i_{2r+1} 
&(0,i_{s}i_{r})
&(i_{r}i_{s},0)\\\midrule 
\end{array}$

\end{enumerate}

\newpage\subsection{$P_{11}   :(a,b)(c,d)=(ca-d^*b,a^*d+cb)$}
       
   For all $r,s.$\\
   
$\begin{array}{lll} P_{11}&i_{2s}&i_{2s+1}\\\toprule i_{2r}&(i_si_r,0)&(0,i_r^*i_s)  \\\midrule i_{2r+1}  
&(0,i_si_r)  
&(-i_s^*i_r,0)\\\midrule \end{array}$

\begin{enumerate}
  \item For $r=s=0$ (The corner)

$\begin{array}{lll} P_{11}&i_{0}&i_{1}\\\toprule i_{0}&1&i_1\\\midrule i_{1}  
&i_1
&-1\\\midrule 
\end{array}$

  \item For $r>s=0$ (The left edge)

$\begin{array}{lll} P_{11}&i_{0}&i_{1}\\\toprule i_{2r}&i_{2r}&-i_{2r+1}\\\midrule i_{2r+1}  
&i_{2r+1}
&-i_{2r}\\\midrule 
\end{array}$

  \item For $s>r=0$ (The top edge)

$\begin{array}{lll} P_{11}&i_{2s}&i_{2s+1}\\\toprule i_{0}&i_{2s}&i_{2s+1}\\\midrule i_{1}  
&i_{2s+1}
&i_{2s}\\\midrule 
\end{array}$

  \item For $r=s>0$ (The diagonal)

$\begin{array}{lll} P_{11}&i_{2r}&i_{2r+1}\\\toprule i_{2r}&-1&i_1\\\midrule i_{2r+1}  
&-i_1
&-1\\\midrule 
\end{array}$
  
  \item For $0\ne r\ne s\ne0$ (The interior)
  
$\begin{array}{lll} P_{11}&i_{2s}&i_{2s+1}\\\toprule i_{2r}&(i_{s}i_{r},0)&-(0,i_{r}i_{s})\\\midrule
i_{2r+1} 
&(0,i_{s}i_{r})
&(i_{s}i_{r},0)\\\midrule 
\end{array}$

\end{enumerate}

\newpage\subsection{$P_{12}  :(a,b)(c,d)=(ac-bd^*,da^*+bc)$}  
For all $r,s.$\\

$\begin{array}{lll} P_{12}&i_{2s}&i_{2s+1}\\\toprule i_{2r}&(i_ri_s,0)&(0,i_si_r^*)  \\\midrule i_{2r+1}  
&(0,i_ri_s)  
&(-i_ri_s^*,0)\\\midrule \end{array}$

\begin{enumerate}
  \item For $r=s=0$ (The corner)

$\begin{array}{lll} P_{12}&i_{0}&i_{1}\\\toprule i_{0}&1&i_1\\\midrule i_{1}  
&i_1
&-1\\\midrule 
\end{array}$

  \item For $r>s=0$ (The left edge)

$\begin{array}{lll} P_{12}&i_{0}&i_{1}\\\toprule i_{2r}&i_{2r}&-i_{2r+1}\\\midrule i_{2r+1}  
&i_{2r+1}
&-i_{2r}\\\midrule 
\end{array}$

  \item For $s>r=0$ (The top edge)

$\begin{array}{lll} P_{12}&i_{2s}&i_{2s+1}\\\toprule i_{0}&i_{2s}&i_{2s+1}\\\midrule i_{1}  
&i_{2s+1}
&i_{2s}\\\midrule 
\end{array}$

  \item For $r=s>0$ (The diagonal)

$\begin{array}{lll} P_{12}&i_{2r}&i_{2r+1}\\\toprule i_{2r}&-1&i_1\\\midrule i_{2r+1}  
&-i_1
&-1\\\midrule 
\end{array}$
  
  \item For $0\ne r\ne s\ne0$ (The interior)
  
$\begin{array}{lll} P_{12}&i_{2s}&i_{2s+1}\\\toprule i_{2r}&(i_{r}i_{s},0)&-(0,i_{s}i_{r})\\\midrule
i_{2r+1} 
&(0,i_{r}i_{s})
&(i_{r}i_{s},0)\\\midrule 
\end{array}$

\end{enumerate}

\newpage\subsection{$P_{13}  :(a,b)(c,d)=(ac-d^*b,da^*+bc)$}  
For all $r,s.$\\

$\begin{array}{lll} P_{13}&i_{2s}&i_{2s+1}\\\toprule i_{2r}&(i_ri_s,0)&(0,i_si_r^*)  \\\midrule i_{2r+1}  
&(0,i_ri_s)  
&(-i_s^*i_r,0)\\\midrule \end{array}$

\begin{enumerate}
  \item For $r=s=0$ (The corner)

$\begin{array}{lll} P_{13}&i_{0}&i_{1}\\\toprule i_{0}&1&i_1\\\midrule i_{1}  
&i_1
&-1\\\midrule 
\end{array}$

  \item For $r>s=0$ (The left edge)

$\begin{array}{lll} P_{13}&i_{0}&i_{1}\\\toprule i_{2r}&i_{2r}&-i_{2r+1}\\\midrule i_{2r+1}  
&i_{2r+1}
&-i_{2r}\\\midrule 
\end{array}$

  \item For $s>r=0$ (The top edge)

$\begin{array}{lll} P_{13}&i_{2s}&i_{2s+1}\\\toprule i_{0}&i_{2s}&i_{2s+1}\\\midrule i_{1}  
&i_{2s+1}
&i_{2s}\\\midrule 
\end{array}$

  \item For $r=s>0$ (The diagonal)

$\begin{array}{lll} P_{13}&i_{2r}&i_{2r+1}\\\toprule i_{2r}&-1&i_1\\\midrule i_{2r+1}  
&-i_1
&-1\\\midrule 
\end{array}$
  
  \item For $0\ne r\ne s\ne0$ (The interior)
  
$\begin{array}{lll} P_{13}&i_{2s}&i_{2s+1}\\\toprule i_{2r}&(i_{r}i_{s},0)&-(0,i_{s}i_{r})\\\midrule
i_{2r+1} 
&(0,i_{r}i_{s})
&(i_{s}i_{r},0)\\\midrule 
\end{array}$

\end{enumerate}

\newpage\subsection{$P_{14}  :(a,b)(c,d)=(ac-bd^*,a^*d+cb)$}  
For all $r,s.$\\

$\begin{array}{lll} P_{14}&i_{2s}&i_{2s+1}\\\toprule i_{2r}&(i_ri_s,0)&(0,i_r^*i_s)  \\\midrule i_{2r+1}  
&(0,i_si_r)  
&(-i_ri_s^*,0)\\\midrule \end{array}$

\begin{enumerate}
  \item For $r=s=0$ (The corner)

$\begin{array}{lll} P_{14}&i_{0}&i_{1}\\\toprule i_{0}&1&i_1\\\midrule i_{1}  
&i_1
&-1\\\midrule 
\end{array}$

  \item For $r>s=0$ (The left edge)

$\begin{array}{lll} P_{14}&i_{0}&i_{1}\\\toprule i_{2r}&i_{2r}&-i_{2r+1}\\\midrule i_{2r+1}  
&i_{2r+1}
&-i_{2r}\\\midrule 
\end{array}$

  \item For $s>r=0$ (The top edge)

$\begin{array}{lll} P_{14}&i_{2s}&i_{2s+1}\\\toprule i_{0}&i_{2s}&i_{2s+1}\\\midrule i_{1}  
&i_{2s+1}
&i_{2s}\\\midrule 
\end{array}$

  \item For $r=s>0$ (The diagonal)

$\begin{array}{lll} P_{14}&i_{2r}&i_{2r+1}\\\toprule i_{2r}&-1&i_1\\\midrule i_{2r+1}  
&-i_1
&-1\\\midrule 
\end{array}$
  
  \item For $0\ne r\ne s\ne0$ (The interior)
  
$\begin{array}{lll} P_{14}&i_{2s}&i_{2s+1}\\\toprule i_{2r}&(i_{r}i_{s},0)&-(0,i_{r}i_{s})\\\midrule
i_{2r+1} 
&(0,i_{s}i_{r})
&(i_{r}i_{s},0)\\\midrule 
\end{array}$

\end{enumerate}

\newpage\subsection{$P_{15}  :(a,b)(c,d)=(ac-d^*b,a^*d+cb)$}  
For all $r,s.$\\

$\begin{array}{lll} P_{15}&i_{2s}&i_{2s+1}\\\toprule i_{2r}&(i_ri_s,0)&(0,i_r^*i_s)  \\\midrule i_{2r+1} 
&(0,i_si_r) 
&(-i_s^*i_r,0)\\\midrule \end{array}$

\begin{enumerate}
  \item For $r=s=0$ (The corner)

$\begin{array}{lll} P_{15}&i_{0}&i_{1}\\\toprule i_{0}&1&i_1\\\midrule i_{1}  
&i_1
&-1\\\midrule 
\end{array}$

  \item For $r>s=0$ (The left edge)

$\begin{array}{lll} P_{15}&i_{0}&i_{1}\\\toprule i_{2r}&i_{2r}&-i_{2r+1}\\\midrule i_{2r+1}  
&i_{2r+1}
&-i_{2r}\\\midrule 
\end{array}$

  \item For $s>r=0$ (The top edge)

$\begin{array}{lll} P_{15}&i_{2s}&i_{2s+1}\\\toprule i_{0}&i_{2s}&i_{2s+1}\\\midrule i_{1}  
&i_{2s+1}
&i_{2s}\\\midrule 
\end{array}$

  \item For $r=s>0$ (The diagonal)

$\begin{array}{lll} P_{15}&i_{2r}&i_{2r+1}\\\toprule i_{2r}&-1&i_1\\\midrule i_{2r+1}  
&-i_1
&-1\\\midrule 
\end{array}$
  
  \item For $0\ne r\ne s\ne0$ (The interior)
  
$\begin{array}{lll} P_{15}&i_{2s}&i_{2s+1}\\\toprule i_{2r}&(i_{r}i_{s},0)&-(0,i_{r}i_{s})\\\midrule
i_{2r+1} 
&(0,i_{s}i_{r})
&(i_{s}i_{r},0)\\\midrule 
\end{array}$

\end{enumerate}

     \newpage\subsection{$P_{16}  :(a,b)(c,d)=(ca-b^*d,ad+c^*b)$}  
For all $r,s.$\\

$\begin{array}{lll} P_{16}&i_{2s}&i_{2s+1}\\\toprule i_{2r}&(i_si_r,0)&(0,i_ri_s)  \\\midrule i_{2r+1}  
&(0,i_s^*i_r)  
&(-i_r^*i_s,0)\\\midrule \end{array}$

\begin{enumerate}
  \item For $r=s=0$ (The corner)

$\begin{array}{lll} P_{16}&i_{0}&i_{1}\\\toprule i_{0}&1&i_1\\\midrule i_{1}  
&i_1
&-1\\\midrule 
\end{array}$

  \item For $r>s=0$ (The left edge)

$\begin{array}{lll} P_{16}&i_{0}&i_{1}\\\toprule i_{2r}&i_{2r}&i_{2r+1}\\\midrule i_{2r+1}  
&i_{2r+1}
&i_{2r}\\\midrule 
\end{array}$

  \item For $s>r=0$ (The top edge)

$\begin{array}{lll} P_{16}&i_{2s}&i_{2s+1}\\\toprule i_{0}&i_{2s}&i_{2s+1}\\\midrule i_{1}  
&-i_{2s+1}
&-i_{2s}\\\midrule 
\end{array}$

  \item For $r=s>0$ (The diagonal)

$\begin{array}{lll} P_{16}&i_{2r}&i_{2r+1}\\\toprule i_{2r}&-1&-i_1\\\midrule i_{2r+1}  
&i_1
&-1\\\midrule 
\end{array}$
  
  \item For $0\ne r\ne s\ne0$ (The interior)
  
$\begin{array}{lll} P_{16}&i_{2s}&i_{2s+1}\\\toprule i_{2r}&(i_{s}i_{r},0)&(0,i_{r}i_{s})\\\midrule
i_{2r+1} 
&-(0,i_{s}i_{r})
&(i_{r}i_{s},0)\\\midrule 
\end{array}$

\end{enumerate}

     \newpage\subsection{$P_{17}   :(a,b)(c,d)=(ca-db^*,ad+c^*b)$}  
For all $r,s.$\\

$\begin{array}{lll} P_{17}&i_{2s}&i_{2s+1}\\\toprule i_{2r}&(i_si_r,0)&(0,i_ri_s)  \\\midrule i_{2r+1}  
&(0,i_s^*i_r)  
&(-i_si_r^*,0)\\\midrule \end{array}$
  
\begin{enumerate}
  \item For $r=s=0$ (The corner)

$\begin{array}{lll} P_{17}&i_{0}&i_{1}\\\toprule i_{0}&1&i_1\\\midrule i_{1}  
&i_1
&-1\\\midrule 
\end{array}$

  \item For $r>s=0$ (The left edge)

$\begin{array}{lll} P_{17}&i_{0}&i_{1}\\\toprule i_{2r}&i_{2r}&i_{2r+1}\\\midrule i_{2r+1}  
&i_{2r+1}
&i_{2r}\\\midrule 
\end{array}$

  \item For $s>r=0$ (The top edge)

$\begin{array}{lll} P_{17}&i_{2s}&i_{2s+1}\\\toprule i_{0}&i_{2s}&i_{2s+1}\\\midrule i_{1}  
&-i_{2s+1}
&-i_{2s}\\\midrule 
\end{array}$

  \item For $r=s>0$ (The diagonal)

$\begin{array}{lll} P_{17}&i_{2r}&i_{2r+1}\\\toprule i_{2r}&-1&-i_1\\\midrule i_{2r+1}  
&i_1
&-1\\\midrule 
\end{array}$
  
  \item For $0\ne r\ne s\ne0$ (The interior)
  
$\begin{array}{lll} P_{17}&i_{2s}&i_{2s+1}\\\toprule i_{2r}&(i_{s}i_{r},0)&(0,i_{r}i_{s})\\\midrule i_{2r+1}
 
&-(0,i_{s}i_{r})
&(i_{s}i_{r},0)\\\midrule 
\end{array}$

\end{enumerate}

     \newpage\subsection{$P_{18}  :(a,b)(c,d)=(ca-b^*d,da+bc^*)$}  
For all $r,s.$\\

$\begin{array}{lll} P_{18}&i_{2s}&i_{2s+1}\\\toprule i_{2r}&(i_si_r,0)&(0,i_si_r)  \\\midrule i_{2r+1}  
&(0,i_ri_s^*)  
d&(-i_r^*i_s,0)\\\midrule \end{array}$

\begin{enumerate}
  \item For $r=s=0$ (The corner)

$\begin{array}{lll} P_{18}&i_{0}&i_{1}\\\toprule i_{0}&1&i_1\\\midrule i_{1}  
&i_1
&-1\\\midrule 
\end{array}$

  \item For $r>s=0$ (The left edge)

$\begin{array}{lll} P_{18}&i_{0}&i_{1}\\\toprule i_{2r}&i_{2r}&i_{2r+1}\\\midrule i_{2r+1}  
&i_{2r+1}
&i_{2r}\\\midrule 
\end{array}$

  \item For $s>r=0$ (The top edge)

$\begin{array}{lll} P_{18}&i_{2s}&i_{2s+1}\\\toprule i_{0}&i_{2s}&i_{2s+1}\\\midrule i_{1}  
&-i_{2s+1}
&-i_{2s}\\\midrule 
\end{array}$

  \item For $r=s>0$ (The diagonal)

$\begin{array}{lll} P_{18}&i_{2r}&i_{2r+1}\\\toprule i_{2r}&-1&-i_1\\\midrule i_{2r+1}  
&i_1
&-1\\\midrule 
\end{array}$
  
  \item For $0\ne r\ne s\ne0$ (The interior)
  
$\begin{array}{lll} P_{18}&i_{2s}&i_{2s+1}\\\toprule i_{2r}&(i_{s}i_{r},0)&(0,i_{s}i_{r})\\\midrule
i_{2r+1} 
&-(0,i_{r}i_{s})
&(i_{r}i_{s},0)\\\midrule 
\end{array}$

\end{enumerate}

     \newpage\subsection{$P_{19}   :(a,b)(c,d)=(ca-db^*,da+bc^*)$}  
For all $r,s.$\\

$\begin{array}{lll} P_{19}&i_{2s}&i_{2s+1}\\\toprule i_{2r}&(i_si_r,0)&(0,i_si_r)  \\\midrule i_{2r+1}  
&(0,i_ri_s^*)  
&(-i_si_r^*,0)\\\midrule \end{array}$
  
\begin{enumerate}
  \item For $r=s=0$ (The corner)

$\begin{array}{lll} P_{19}&i_{0}&i_{1}\\\toprule i_{0}&1&i_1\\\midrule i_{1}  
&i_1
&-1\\\midrule 
\end{array}$

  \item For $r>s=0$ (The left edge)

$\begin{array}{lll} P_{19}&i_{0}&i_{1}\\\toprule i_{2r}&i_{2r}&i_{2r+1}\\\midrule i_{2r+1}  
&i_{2r+1}
&i_{2r}\\\midrule 
\end{array}$

  \item For $s>r=0$ (The top edge)

$\begin{array}{lll} P_{19}&i_{2s}&i_{2s+1}\\\toprule i_{0}&i_{2s}&i_{2s+1}\\\midrule i_{1}  
&-i_{2s+1}
&-i_{2s}\\\midrule 
\end{array}$

  \item For $r=s>0$ (The diagonal)

$\begin{array}{lll} P_{19}&i_{2r}&i_{2r+1}\\\toprule i_{2r}&-1&-i_1\\\midrule i_{2r+1}  
&i_1
&-1\\\midrule 
\end{array}$
  
  \item For $0\ne r\ne s\ne0$ (The interior)
  
$\begin{array}{lll} P_{19}&i_{2s}&i_{2s+1}\\\toprule i_{2r}&(i_{s}i_{r},0)&(0,i_{s}i_{r})\\\midrule i_{2r+1}
 
&-(0,i_{r}i_{s})
&(i_{s}i_{r},0)\\\midrule 
\end{array}$

\end{enumerate}

\newpage\subsection{$P_{20}  :(a,b)(c,d)=(ac-b^*d,ad+c^*b)$}  
For all $r,s.$\\

$\begin{array}{lll} P_{20}&i_{2s}&i_{2s+1}\\\toprule i_{2r}&(i_ri_s,0)&(0,i_ri_s)  \\\midrule i_{2r+1}  
&(0,i_s^*i_r)  
&(-i_r^*i_s,0)\\\midrule \end{array}$

\begin{enumerate}
  \item For $r=s=0$ (The corner)

$\begin{array}{lll} P_{20}&i_{0}&i_{1}\\\toprule i_{0}&1&i_1\\\midrule i_{1}  
&i_1
&-1\\\midrule 
\end{array}$

  \item For $r>s=0$ (The left edge)

$\begin{array}{lll} P_{20}&i_{0}&i_{1}\\\toprule i_{2r}&i_{2r}&i_{2r+1}\\\midrule i_{2r+1}  
&i_{2r+1}
&i_{2r}\\\midrule 
\end{array}$

  \item For $s>r=0$ (The top edge)

$\begin{array}{lll} P_{20}&i_{2s}&i_{2s+1}\\\toprule i_{0}&i_{2s}&i_{2s+1}\\\midrule i_{1}  
&-i_{2s+1}
&-i_{2s}\\\midrule 
\end{array}$

  \item For $r=s>0$ (The diagonal)

$\begin{array}{lll} P_{20}&i_{2r}&i_{2r+1}\\\toprule i_{2r}&-1&-i_1\\\midrule i_{2r+1}  
&i_1
&-1\\\midrule 
\end{array}$
  
  \item For $0\ne r\ne s\ne0$ (The interior)
  
$\begin{array}{lll} P_{20}&i_{2s}&i_{2s+1}\\\toprule i_{2r}&(i_{r}i_{s},0)&(0,i_{r}i_{s})\\\midrule
i_{2r+1} 
&-(0,i_{s}i_{r})
&(i_{r}i_{s},0)\\\midrule 
\end{array}$

\end{enumerate}

\newpage\subsection{$P_{21}  :(a,b)(c,d)=(ac-db^*,ad+c^*b)$}  
For all $r,s.$\\

$\begin{array}{lll} P_{21}&i_{2s}&i_{2s+1}\\\toprule i_{2r}&(i_ri_s,0)&(0,i_ri_s)  \\\midrule i_{2r+1}  
&(0,i_s^*i_r)  
&(-i_si_r^*,0)\\\midrule \end{array}$

\begin{enumerate}
  \item For $r=s=0$ (The corner)

$\begin{array}{lll} P_{21}&i_{0}&i_{1}\\\toprule i_{0}&1&i_1\\\midrule i_{1}  
&i_1
&-1\\\midrule 
\end{array}$

  \item For $r>s=0$ (The left edge)

$\begin{array}{lll} P_{21}&i_{0}&i_{1}\\\toprule i_{2r}&i_{2r}&i_{2r+1}\\\midrule i_{2r+1}  
&i_{2r+1}
&i_{2r}\\\midrule 
\end{array}$

  \item For $s>r=0$ (The top edge)

$\begin{array}{lll} P_{21}&i_{2s}&i_{2s+1}\\\toprule i_{0}&i_{2s}&i_{2s+1}\\\midrule i_{1}  
&-i_{2s+1}
&-i_{2s}\\\midrule 
\end{array}$

  \item For $r=s>0$ (The diagonal)

$\begin{array}{lll} P_{21}&i_{2r}&i_{2r+1}\\\toprule i_{2r}&-1&-i_1\\\midrule i_{2r+1}  
&i_1
&-1\\\midrule 
\end{array}$
  
  \item For $0\ne r\ne s\ne0$ (The interior)
  
$\begin{array}{lll} P_{21}&i_{2s}&i_{2s+1}\\\toprule i_{2r}&(i_{r}i_{s},0)&(0,i_{r}i_{s})\\\midrule
i_{2r+1} 
&-(0,i_{s}i_{r})
&(i_{s}i_{r},0)\\\midrule 
\end{array}$

\end{enumerate}

     \newpage\subsection{$P_{22}  :(a,b)(c,d)=(ac-b^*d,da+bc^*)$}  
For all $r,s.$\\

$\begin{array}{lll} P_{22}&i_{2s}&i_{2s+1}\\\toprule i_{2r}&(i_ri_s,0)&(0,i_si_r)  \\\midrule i_{2r+1}  
&(0,i_ri_s^*)  
&(-i_r^*i_s,0)\\\midrule \end{array}$

\begin{enumerate}
  \item For $r=s=0$ (The corner)

$\begin{array}{lll} P_{22}&i_{0}&i_{1}\\\toprule i_{0}&1&i_1\\\midrule i_{1}  
&i_1
&-1\\\midrule 
\end{array}$

  \item For $r>s=0$ (The left edge)

$\begin{array}{lll} P_{22}&i_{0}&i_{1}\\\toprule i_{2r}&i_{2r}&i_{2r+1}\\\midrule i_{2r+1}  
&i_{2r+1}
&i_{2r}\\\midrule 
\end{array}$

  \item For $s>r=0$ (The top edge)

$\begin{array}{lll} P_{22}&i_{2s}&i_{2s+1}\\\toprule i_{0}&i_{2s}&i_{2s+1}\\\midrule i_{1}  
&-i_{2s+1}
&-i_{2s}\\\midrule 
\end{array}$

  \item For $r=s>0$ (The diagonal)

$\begin{array}{lll} P_{22}&i_{2r}&i_{2r+1}\\\toprule i_{2r}&-1&-i_1\\\midrule i_{2r+1}  
&i_1
&-1\\\midrule 
\end{array}$
  
  \item For $0\ne r\ne s\ne0$ (The interior)
  
$\begin{array}{lll} P_{22}&i_{2s}&i_{2s+1}\\\toprule i_{2r}&(i_{r}i_{s},0)&(0,i_{s}i_{r})\\\midrule
i_{2r+1} 
&-(0,i_{r}i_{s})
&(i_{r}i_{s},0)\\\midrule 
\end{array}$

\end{enumerate}

     \newpage\subsection{$P_{23}  :(a,b)(c,d)=(ac-db^*,da+bc^*)$}  
For all $r,s.$\\

$\begin{array}{lll} P_{23}&i_{2s}&i_{2s+1}\\\toprule i_{2r}&(i_ri_s,0)&(0,i_si_r)  \\\midrule i_{2r+1}  
&(0,i_ri_s^*)  
&(-i_si_r^*,0)\\\midrule \end{array}$

\begin{enumerate}
  \item For $r=s=0$ (The corner)

$\begin{array}{lll} P_{23}&i_{0}&i_{1}\\\toprule i_{0}&1&i_1\\\midrule i_{1}  
&i_1
&-1\\\midrule 
\end{array}$

  \item For $r>s=0$ (The left edge)

$\begin{array}{lll} P_{23}&i_{0}&i_{1}\\\toprule i_{2r}&i_{2r}&i_{2r+1}\\\midrule i_{2r+1}  
&i_{2r+1}
&i_{2r}\\\midrule 
\end{array}$

  \item For $s>r=0$ (The top edge)

$\begin{array}{lll} P_{23}&i_{2s}&i_{2s+1}\\\toprule i_{0}&i_{2s}&i_{2s+1}\\\midrule i_{1}  
&-i_{2s+1}
&-i_{2s}\\\midrule 
\end{array}$

  \item For $r=s>0$ (The diagonal)

$\begin{array}{lll} P_{23}&i_{2r}&i_{2r+1}\\\toprule i_{2r}&-1&i_1\\\midrule i_{2r+1}  
&-i_1
&-1\\\midrule 
\end{array}$
  
  \item For $0\ne r\ne s\ne0$ (The interior)
  
$\begin{array}{lll} P_{23}&i_{2s}&i_{2s+1}\\\toprule i_{2r}&(i_{r}i_{s},0)&(0,i_{s}i_{r})\\\midrule
i_{2r+1} 
&-(0,i_{r}i_{s})
&(i_{s}i_{r},0)\\\midrule 
\end{array}$

\end{enumerate}

     \newpage\subsection{$P_{24}   :(a,b)(c,d)=(ca-bd^*,ad+c^*b)$}  
For all $r,s.$\\

$\begin{array}{lll} P_{24}&i_{2s}&i_{2s+1}\\\toprule i_{2r}&(i_si_r,0)&(0,i_ri_s)  \\\midrule i_{2r+1}  
&(0,i_s^*i_r)  
&(-i_ri_s^*,0)\\\midrule \end{array}$
  
\begin{enumerate}
  \item For $r=s=0$ (The corner)

$\begin{array}{lll} P_{24}&i_{0}&i_{1}\\\toprule i_{0}&1&i_1\\\midrule i_{1}  
&i_1
&-1\\\midrule 
\end{array}$

  \item For $r>s=0$ (The left edge)

$\begin{array}{lll} P_{24}&i_{0}&i_{1}\\\toprule i_{2r}&i_{2r}&i_{2r+1}\\\midrule i_{2r+1}  
&i_{2r+1}
&-i_{2r}\\\midrule 
\end{array}$

  \item For $s>r=0$ (The top edge)

$\begin{array}{lll} P_{24}&i_{2s}&i_{2s+1}\\\toprule i_{0}&i_{2s}&i_{2s+1}\\\midrule i_{1}  
&-i_{2s+1}
&i_{2s}\\\midrule 
\end{array}$

  \item For $r=s>0$ (The diagonal)

$\begin{array}{lll} P_{24}&i_{2r}&i_{2r+1}\\\toprule i_{2r}&-1&-i_1\\\midrule i_{2r+1}  
&i_1
&-1\\\midrule 
\end{array}$
  
  \item For $0\ne r\ne s\ne0$ (The interior)
  
$\begin{array}{lll} P_{24}&i_{2s}&i_{2s+1}\\\toprule i_{2r}&(i_{s}i_{r},0)&(0,i_{r}i_{s})\\\midrule i_{2r+1}
 
&-(0,i_{s}i_{r})
&(i_{r}i_{s},0)\\\midrule 
\end{array}$

\end{enumerate}

     \newpage\subsection{$P_{25}   :(a,b)(c,d)=(ca-d^*b,ad+c^*b)$}
       
For all $r,s.$\\

$\begin{array}{lll} P_{25}&i_{2s}&i_{2s+1}\\\toprule i_{2r}&(i_si_r,0)&(0,i_ri_s)  \\\midrule i_{2r+1}  
&(0,i_s^*i_r)  
&(-i_s^*i_r,0)\\\midrule \end{array}$

\begin{enumerate}
  
  \item For $r=s=0$ (The corner)

$\begin{array}{lll} P_{25}&i_{0}&i_{1}\\\toprule i_{0}&1&i_1\\\midrule i_{1}  
&i_1
&-1\\\midrule 
\end{array}$

  \item For $r>s=0$ (The left edge)

$\begin{array}{lll} P_{25}&i_{0}&i_{1}\\\toprule i_{2r}&i_{2r}&i_{2r+1}\\\midrule i_{2r+1}  
&i_{2r+1}
&-i_{2r}\\\midrule 
\end{array}$

  \item For $s>r=0$ (The top edge)

$\begin{array}{lll} P_{25}&i_{2s}&i_{2s+1}\\\toprule i_{0}&i_{2s}&i_{2s+1}\\\midrule i_{1}  
&-i_{2s+1}
&i_{2s}\\\midrule 
\end{array}$

  \item For $r=s>0$ (The diagonal)

$\begin{array}{lll} P_{25}&i_{2r}&i_{2r+1}\\\toprule i_{2r}&-1&-i_1\\\midrule i_{2r+1}  
&i_1
&-1\\\midrule 
\end{array}$
  
  \item For $0\ne r\ne s\ne0$ (The interior)

$\begin{array}{lll} P_{25}&i_{2s}&i_{2s+1}\\\toprule i_{2r}&(i_{s}i_{r},0)&(0,i_{r}i_{s})\\\midrule i_{2r+1}
 
&-(0,i_{s}i_{r})
&(i_{s}i_{r},0)\\\midrule 
\end{array}$

\end{enumerate}

     \newpage\subsection{$P_{26}   :(a,b)(c,d)=(ca-bd^*,da+bc^*)$}  
For all $r,s.$\\

$\begin{array}{lll} P_{26}&i_{2s}&i_{2s+1}\\\toprule i_{2r}&(i_si_r,0)&(0,i_si_r)  \\\midrule i_{2r+1}  
&(0,i_ri_s^*)  
&(-i_ri_s^*,0)\\\midrule \end{array}$
  
\begin{enumerate}
  \item For $r=s=0$ (The corner)

$\begin{array}{lll} P_{26}&i_{0}&i_{1}\\\toprule i_{0}&1&i_1\\\midrule i_{1}  
&i_1
&-1\\\midrule 
\end{array}$

  \item For $r>s=0$ (The left edge)

$\begin{array}{lll} P_{26}&i_{0}&i_{1}\\\toprule i_{2r}&i_{2r}&i_{2r+1}\\\midrule i_{2r+1}  
&i_{2r+1}
&-i_{2r}\\\midrule 
\end{array}$

  \item For $s>r=0$ (The top edge)

$\begin{array}{lll} P_{26}&i_{2s}&i_{2s+1}\\\toprule i_{0}&i_{2s}&i_{2s+1}\\\midrule i_{1}  
&-i_{2s+1}
&i_{2s}\\\midrule 
\end{array}$

  \item For $r=s>0$ (The diagonal)

$\begin{array}{lll} P_{26}&i_{2r}&i_{2r+1}\\\toprule i_{2r}&-1&-i_1\\\midrule i_{2r+1}  
&i_1
&-1\\\midrule 
\end{array}$
  
  \item For $0\ne r\ne s\ne0$ (The interior)
  
$\begin{array}{lll} P_{26}&i_{2s}&i_{2s+1}\\\toprule i_{2r}&(i_{s}i_{r},0)&(0,i_{s}i_{r})\\\midrule i_{2r+1}
 
&-(0,i_{r}i_{s})
&(i_{r}i_{s},0)\\\midrule 
\end{array}$

\end{enumerate}

\newpage\subsection{$P_{27}   :(a,b)(c,d)=(ca-d^*b,da+bc^*)$}
For all $r,s.$\\

$\begin{array}{lll} P_{27}&i_{2s}&i_{2s+1}\\\toprule i_{2r}&(i_si_r,0)&(0,i_si_r)\\\midrule i_{2r+1}  
&(0,i_ri_s^*)
&(-i_s^*i_r,0)\\\midrule 
\end{array}$
     \begin{enumerate}
       \item For $r=s=0$ (The corner)

$\begin{array}{lll} P_{27}&i_{0}&i_{1}\\\toprule i_{0}&1&i_1\\\midrule i_{1}  
&i_1
&-1\\\midrule 
\end{array}$

       \item For $r>s=0$ (The left edge)

$\begin{array}{lll} P_{27}&i_{0}&i_{1}\\\toprule i_{2r}&i_{2r}&i_{2r+1}\\\midrule i_{2r+1}  
&i_{2r+1}
&-i_{2r}\\\midrule 
\end{array}$

       \item For $s>r=0$ (The top edge)

$\begin{array}{lll} P_{27}&i_{2s}&i_{2s+1}\\\toprule i_{0}&i_{2s}&i_{2s+1}\\\midrule i_{1}  
&-i_{2s+1}
&i_{2s}\\\midrule 
\end{array}$

       \item For $r=s>0$ (The diagonal)

$\begin{array}{lll} P_{27}&i_{2r}&i_{2r+1}\\\toprule i_{2r}&-1&-i_1\\\midrule i_{2r+1}  
&i_1
&-1\\\midrule 
\end{array}$

        \item For $0\ne r\ne s\ne0$ (The interior)

$\begin{array}{lll} P_{27}&i_{2s}&i_{2s+1}\\\toprule i_{2r}&(i_si_r,0)&(0,i_si_r)\\\midrule i_{2r+1}  
&-(0,i_ri_s)
&(i_si_r,0)\\\midrule 
\end{array}$

\end{enumerate}
 
\newpage\subsection{$P_{28}  :(a,b)(c,d)=(ac-bd^*,ad+c^*b)$}  
For all $r,s.$\\

$\begin{array}{lll} P_{28}&i_{2s}&i_{2s+1}\\\toprule i_{2r}&(i_ri_s,0)&(0,i_ri_s)  \\\midrule i_{2r+1}  
&(0,i_s^*i_r)  
&(-i_ri_s^*,0)\\\midrule \end{array}$

\begin{enumerate}
  \item For $r=s=0$ (The corner)

$\begin{array}{lll} P_{28}&i_{0}&i_{1}\\\toprule i_{0}&1&i_1\\\midrule i_{1}  
&i_1
&-1\\\midrule 
\end{array}$

  \item For $r>s=0$ (The left edge)

$\begin{array}{lll} P_{28}&i_{0}&i_{1}\\\toprule i_{2r}&i_{2r}&i_{2r+1}\\\midrule i_{2r+1}  
&i_{2r+1}
&-i_{2r}\\\midrule 
\end{array}$

  \item For $s>r=0$ (The top edge)

$\begin{array}{lll} P_{28}&i_{2s}&i_{2s+1}\\\toprule i_{0}&i_{2s}&i_{2s+1}\\\midrule i_{1}  
&-i_{2s+1}
&i_{2s}\\\midrule 
\end{array}$

  \item For $r=s>0$ (The diagonal)

$\begin{array}{lll} P_{28}&i_{2r}&i_{2r+1}\\\toprule i_{2r}&-1&-i_1\\\midrule i_{2r+1}  
&i_1
&-1\\\midrule 
\end{array}$
  
  \item For $0\ne r\ne s\ne0$ (The interior)
  
$\begin{array}{lll} P_{28}&i_{2s}&i_{2s+1}\\\toprule i_{2r}&(i_{r}i_{s},0)&(0,i_{r}i_{s})\\\midrule
i_{2r+1} 
&-(0,i_{s}i_{r})
&(i_{r}i_{s},0)\\\midrule 
\end{array}$

\end{enumerate}

\newpage\subsection{$P_{29}  :(a,b)(c,d)=(ac-d^*b,ad+c^*b)$}  
For all $r,s.$\\

$\begin{array}{lll} P_{29}&i_{2s}&i_{2s+1}\\\toprule i_{2r}&(i_ri_s,0)&(0,i_ri_s)  \\\midrule i_{2r+1}  
&(0,i_s^*i_r)  
&(-i_s^*i_r,0)\\\midrule \end{array}$

\begin{enumerate}
  \item For $r=s=0$ (The corner)

$\begin{array}{lll} P_{29}&i_{0}&i_{1}\\\toprule i_{0}&1&i_1\\\midrule i_{1}  
&i_1
&-1\\\midrule 
\end{array}$

  \item For $r>s=0$ (The left edge)

$\begin{array}{lll} P_{29}&i_{0}&i_{1}\\\toprule i_{2r}&i_{2r}&-i_{2r+1}\\\midrule i_{2r+1}  
&i_{2r+1}
&i_{2r}\\\midrule 
\end{array}$

  \item For $s>r=0$ (The top edge)

$\begin{array}{lll} P_{29}&i_{2s}&i_{2s+1}\\\toprule i_{0}&i_{2s}&i_{2s+1}\\\midrule i_{1}  
&-i_{2s+1}
&i_{2s}\\\midrule 
\end{array}$

  \item For $r=s>0$ (The diagonal)

$\begin{array}{lll} P_{29}&i_{2r}&i_{2r+1}\\\toprule i_{2r}&-1&-i_1\\\midrule i_{2r+1}  
&i_1
&-1\\\midrule 
\end{array}$
  
  \item For $0\ne r\ne s\ne0$ (The interior)
  
$\begin{array}{lll} P_{29}&i_{2s}&i_{2s+1}\\\toprule i_{2r}&(i_{r}i_{s},0)&(0,i_{r}i_{s})\\\midrule i_{2r+1}
 
&-(0,i_{s}i_{r})
&(i_{s}i_{r},0)\\\midrule 
\end{array}$

\end{enumerate}

\newpage\subsection{$P_{30}  :(a,b)(c,d)=(ac-bd^*,da+bc^*)$}  
For all $r,s.$\\

$\begin{array}{lll} P_{30}&i_{2s}&i_{2s+1}\\\toprule i_{2r}&(i_ri_s,0)&(0,i_si_r)  \\\midrule i_{2r+1}  
&(0,i_ri_s^*)  
&(-i_ri_s^*,0)\\\midrule \end{array}$

\begin{enumerate}
  \item For $r=s=0$ (The corner)

$\begin{array}{lll} P_{30}&i_{0}&i_{1}\\\toprule i_{0}&1&i_1\\\midrule i_{1}  
&i_1
&-1\\\midrule 
\end{array}$

  \item For $r>s=0$ (The left edge)

$\begin{array}{lll} P_{30}&i_{0}&i_{1}\\\toprule i_{2r}&i_{2r}&i_{2r+1}\\\midrule i_{2r+1}  
&i_{2r+1}
&-i_{2r}\\\midrule 
\end{array}$

  \item For $s>r=0$ (The top edge)

$\begin{array}{lll} P_{30}&i_{2s}&i_{2s+1}\\\toprule i_{0}&i_{2s}&i_{2s+1}\\\midrule i_{1}  
&-i_{2s+1}
&i_{2s}\\\midrule 
\end{array}$

  \item For $r=s>0$ (The diagonal)

$\begin{array}{lll} P_{9}{0}&i_{2r}&i_{2r+1}\\\toprule i_{2r}&-1&-i_1\\\midrule i_{2r+1}  
&i_1
&-1\\\midrule 
\end{array}$
  
  \item For $0\ne r\ne s\ne0$ (The interior)
  
$\begin{array}{lll} P_{30}&i_{2s}&i_{2s+1}\\\toprule i_{2r}&(i_{r}i_{s},0)&(0,i_{s}i_{r})\\\midrule
i_{2r+1} 
&-(0,i_{r}i_{s})
&(i_{r}i_{s},0)\\\midrule 
\end{array}$

\end{enumerate}

\newpage\subsection{$P_{31}  :(a,b)(c,d)=(ac-d^*b,da+bc^*)$}  
For all $r,s.$\\

$\begin{array}{lll} P_{31}&i_{2s}&i_{2s+1}\\\toprule i_{2r}&(i_ri_s,0)&(0,i_si_r)  \\\midrule i_{2r+1}  
&(0,i_ri_s^*)  
&(-i_s^*i_r,0)\\\midrule \end{array}$

\begin{enumerate}
  \item For $r=s=0$ (The corner)

$\begin{array}{lll} P_{31}&i_{0}&i_{1}\\\toprule i_{0}&1&i_1\\\midrule i_{1}  
&i_1
&-1\\\midrule 
\end{array}$

  \item For $r>s=0$ (The left edge)

$\begin{array}{lll} P_{31}&i_{0}&i_{1}\\\toprule i_{2r}&i_{2r}&-i_{2r+1}\\\midrule i_{2r+1}  
&i_{2r+1}
&i_{2r}\\\midrule 
\end{array}$

  \item For $s>r=0$ (The top edge)

$\begin{array}{lll} P_{31}&i_{2s}&i_{2s+1}\\\toprule i_{0}&i_{2s}&i_{2s+1}\\\midrule i_{1}  
&-i_{2s+1}
&i_{2s}\\\midrule 
\end{array}$

  \item For $r=s>0$ (The diagonal)

$\begin{array}{lll} P_{31}&i_{2r}&i_{2r+1}\\\toprule i_{2r}&-1&-i_1\\\midrule i_{2r+1}  
&i_1
&-1\\\midrule 
\end{array}$
  
  \item For $0\ne r\ne s\ne0$ (The interior)
  
$\begin{array}{lll} P_{31}&i_{2s}&i_{2s+1}\\\toprule i_{2r}&(i_{r}i_{s},0)&(0,i_{s}i_{r})\\\midrule
i_{2r+1} 
&-(0,i_{r}i_{s})
&(i_{s}i_{r},0)\\\midrule 
\end{array}$

\end{enumerate}

\newpage

\section{Cayley-Dickson-like products as twisted group products}

The basis elements $i_0,i_1,i_2,\cdots$ are indexed by $W=\{0,1,2,\cdots\}$ which is a group under the bit-wise
`exclusive or' of their binary representations and for each of the 32 Cayley-Dickson-like products there is a function
$\omega$ from $W\times W$ to $\{-1,1\}$ such that for $p,q\in W$

\begin{equation}
  i_pi_q=\omega(p,q)i_{pq}
\end{equation}

where $pq$ is the group product of $p$ and $q$. The function $\omega$ is called a `twist' on the group $W$ and turns
the set of all finite sequences into a twisted group algebra.

Since $1=(1,0)=i_0$ is the identity, it follows that $\omega(p,0)=\omega(0,q)=1$ for all Cayley-Dickson-like products, and
since for $p>1,$ $i_pi_p^*=\|i_p\|^2=1$ and since $i_pi_p^*=-i_pi_p=-\omega(p,p)i_0=1$ it follows that for $p>1$

\begin{equation}
  \omega(p,p)=-1
\end{equation}

and that for all $p$

\begin{equation}
  i_p^*=\omega(p,p)i_p
\end{equation}

From equation \vref{eq:neg} it follows that

\begin{equation}\label{Eqn:neg}
  \omega(q,p)+\omega(p,q)=0\text{\ for\ }0\ne p\ne q\ne0
\end{equation}

\subsection{The product of $i_{2r}i_{2s}$}

Since $i_{2r}=\left(i_r,0\right)$ then $i_{2r}i_{2s}=\left(i_r,0\right)\left(i_s,0\right)$ this could be called the
case of $b=d=0.$ For each of the 32 distinct products, either $i_{2r}i_{2s}=\left(i_si_r,0\right)$ or
$i_{2r}i_{2s}=\left(i_ri_s,0\right)$. We shall consider the effect each of these alternatives upon the twist $\omega$.

 \subsubsection{$i_{2r}i_{2s}=\left(i_si_r,0\right)$}
 
  Since $i_{2r}i_{2s}=\omega(2r,2s)i_{2rs}$ and since $i_{2r}i_{2s}=\omega(s,r)i_{2rs}$ we may conclude that
Whenever $i_{2r}i_{2s}=\left(i_si_r,0\right)$,
  \begin{equation}
     \omega(2r,2s)=\omega(s,r)
  \end{equation}

 \subsubsection{$i_{2r}i_{2s}=\left(i_ri_s,0\right)$}
 
   Since $i_{2r}i_{2s}=\omega(2r,2s)i_{2rs}$ and since $i_{2r}i_{2s}=\omega(r,s)i_{2rs}$ we may conclude that
Whenever $i_{2r}i_{2s}=\left(i_ri_s,0\right)$,
  \begin{equation}
     \omega(2r,2s)=\omega(r,s)
  \end{equation}

\subsection{The products of $i_{2r}i_{2s+1}$ and $i_{2r+1}i_{2s}$}

Since $i_{2r}i_{2s+1}=\left(i_r,0\right)\left(0,i_s\right)$ this could be called the case of $b=c=0$. And since
$i_{2r+1}i_{2s}=\left(0,i_r\right)\left(i_s,0\right)$ that could be called the case of $a=d=0.$ An inspection of the 32
alternate products shows that there are only four distinct possibilities for the products of  $i_{2r}i_{2s+1}$ and
$i_{2r+1}i_{2s}$.\\

 \subsubsection{ $i_{2r}i_{2s+1 }=\left(0 ,i_si_r  \right)$ and $i_{2r+1 }i_{2s }=\left(0,i_ri_s^* 
\right)$}
 
 These two conditions imply that whenever $i_{2r}i_{2s+1 }=\left(0 ,i_si_r  \right)$,
 
 \begin{equation}
    \omega(2r,2s+1)=\omega(s,r)
 \end{equation}
and that whenever $i_{2r+1 }i_{2s }=\left(i_ri_s^* ,0  \right)$,
\begin{equation}
  \omega(2r+1,2s)=\omega(s,s)\omega(r,s)=\begin{cases}-\omega(r,s)\text{\ if\ }s>0\\
                                                      1\text{\ otherwise}
                                         \end{cases}
\end{equation}

 \subsubsection{$i_{2r}i_{2s+1 }=\left(0 ,i_ri_s  \right)$ and $i_{2r+1 }i_{2s }=\left(0,i_s^*i_r 
\right)$}
 
 These two conditions imply that whenever $i_{2r}i_{2s+1 }=\left(0,i_ri_s  \right)$,
 
 \begin{equation}
    \omega(2r,2s+1)=\omega(r,s)
 \end{equation}
and that whenever $i_{2r+1 }i_{2s }=\left(i_s^*i_r ,0  \right)$,
\begin{equation}
  \omega(2r+1,2s)=\omega(s,s)\omega(s,r)=\begin{cases}-\omega(s,r)\text{\ if\ }s>0\\
                                                      1\text{\ otherwise}
                                         \end{cases}
\end{equation}

 \subsubsection{ $i_{2r}i_{2s+1 }=\left(0 ,i_si_r^*  \right)$ and $i_{2r+1 }i_{2s }=\left(0,i_ri_s 
\right)$}
 
These two conditions imply that whenever $i_{2r}i_{2s+1 }=\left(0 ,i_si_r^*  \right)$
 
 \begin{equation}
    \omega(2r,2s+1)=\omega(r,r)\omega(s,r)=\begin{cases}
                                             -\omega(s,r)\text{\ if\ }r>0\\
                                             1\text{\ otherwise}
                                           \end{cases}
 \end{equation}
 
and that whenever $i_{2r+1 }i_{2s }=\left(i_ri_s ,0  \right)$,

\begin{equation}
  \omega(2r+1,2s)=\omega(r,s)
\end{equation}

 \subsubsection{$i_{2r}i_{2s+1 }=\left(0 ,i_r^*i_s  \right)$ and $i_{2r+1 }i_{2s }=\left(0,i_si_r 
\right)$}
 
 These two conditions imply that whenever $i_{2r}i_{2s+1 }=\left(0 ,i_r^*i_s  \right)$
 
 \begin{equation}
    \omega(2r,2s+1)=\omega(s,s)\omega(r,s)=\begin{cases}
                                             -\omega(r,s)\text{\ if\ }s>0\\
                                             1\text{\ otherwise}
                                           \end{cases}
 \end{equation}
 
and that whenever  $i_{2r+1 }i_{2s }=\left(i_si_r ,0  \right)$,

\begin{equation}
  \omega(2r+1,2s)=\omega(r,s)
\end{equation}

\subsection{The product of $i_{2r+1}i_{2s+1}$}

Since $i_{2r+1}=\left(0,i_r\right)$ the product $i_{2r+1}i_{2s+1}$ could be called the case of $a=c=0.$

 \subsubsection{$i_{2r+1 }i_{2s+1 }=-\left(i_s^*i_r,0  \right)$}
 
 This implies that whenever $i_{2r+1 }i_{2s+1 }=-\left(i_s^*i_r,0  \right)$
 
 \begin{equation}
   \omega(2r+1,2s+1)=-\omega(s,s)\omega(s,r)=\begin{cases}
                                               \omega(s,r)\text{\ if\ }s>0\\
                                             -1\text{\ otherwise}  
                                             \end{cases}
 \end{equation}

 \subsubsection{$i_{2r+1 }i_{2s+1 }=-\left(i_ri_s^*,0  \right)$}
 
 This implies that whenever $i_{2r+1 }i_{2s+1 }=-\left(i_ri_s^*,0)  \right)$
 
  \begin{equation}
   \omega(2r+1,2s+1)=-\omega(s,s)\omega(r,s)=\begin{cases}
                                               \omega(r,s)\text{\ if\ }s>0\\
                                             -1\text{\ otherwise}  
                                             \end{cases}
 \end{equation}

 \subsubsection{$i_{2r+1 }i_{2s+1 }=-\left(i_si_r^*,0  \right)$}
 
 This implies that whenever $i_{2r+1 }i_{2s+1 }=-\left(i_si_r^*,0  \right)$
 
 \begin{equation}
   \omega(2r+1,2s+1)=-\omega(r,r)\omega(s,r)=\begin{cases}
                                               \omega(s,r)\text{\ if\ }r>0\\
                                             -1\text{\ otherwise}  
                                             \end{cases}
 \end{equation}

 \subsubsection{$i_{2r+1 }i_{2s+1 }=-\left(i_r^*i_s,0  \right)$}
 
This implies that whenever $i_{2r+1 }i_{2s+1 }=-\left(i_r^*i_s,0  \right)$

 \begin{equation}
   \omega(2r+1,2s+1)=-\omega(r,r)\omega(r,s)=\begin{cases}
                                               \omega(r,s)\text{\ if\ }r>0\\
                                             -1\text{\ otherwise}  
                                             \end{cases}
 \end{equation}

 \section{The 32 variations on $\omega$}\label{sec:variations}

For each of the 32 Cayley-Dickson-like products of finite sequences, the recursive definitions for $\omega$
are given.

\begin{enumerate}
\setcounter{enumi}{-1}
     \item \fbox{$P_{0}  :(a,b)(c,d)=(ca-b^*d,da^*+bc)$}\\  
\fbox{$\begin{array}{lll}\omega_{0}&2s&2s+1\\\toprule 2r&\omega(s,r)&\begin{cases}-\omega(s,r)\text{\ if\
}r>0\\1\text{\ otherwise}\end{cases}  \\\midrule 2r+1  
&\omega(r,s)  
&\begin{cases}\omega(r,s)\text{\ if\ }r>0\\-1\text{\ otherwise}\end{cases}\end{array}$}
     \item \fbox{$P_{1}  :(a,b)(c,d)=(ca-db^*,da^*+bc)$}\\  
\fbox{$\begin{array}{lll}\omega_{1}&2s&2s+1\\\toprule 2r&\omega(s,r)&\begin{cases}-\omega(s,r)\text{\ if\
}r>0\\1\text{\ otherwise}\end{cases}  \\\midrule 2r+1  
&\omega(r,s)  
&\begin{cases}\omega(s,r)\text{\ if\ }r>0\\-1\text{\ otherwise}\end{cases}\end{array}$}
     \item \fbox{$P_{2}  :(a,b)(c,d)=(ca-b^*d,a^*d+cb)$}\\  
\fbox{$\begin{array}{lll}\omega_{2}&2s&2s+1\\\toprule 2r&\omega(s,r)&\begin{cases}-\omega(r,s)\text{\ if\
}r>0\\1\text{\ otherwise}\end{cases}  \\\midrule 2r+1  
&\omega(s,r)  
&\begin{cases}\omega(r,s)\text{\ if\ }r>0\\-1\text{\ otherwise}\end{cases}\end{array}$}
     \item \fbox{$P_{3}  :(a,b)(c,d)=(ca-db^*,a^*d+cb)$}\\  
\fbox{$\begin{array}{lll}\omega_{3}&2s&2s+1\\\toprule 2r&\omega(s,r)&\begin{cases}-\omega(r,s)\text{\ if\
}r>0\\1\text{\ otherwise}\end{cases}  \\\midrule 2r+1  
&\omega(s,r)  
&\begin{cases}\omega(s,r)\text{\ if\ }r>0\\-1\text{\ otherwise}\end{cases}\end{array}$}
     \item \fbox{$P_{4}  :(a,b)(c,d)=(ac-b^*d,da^*+bc)$}\\  
\fbox{$\begin{array}{lll}\omega_{4}&2s&2s+1\\\toprule 2r&\omega(r,s)&\begin{cases}-\omega(s,r)\text{\ if\
}r>0\\1\text{\ otherwise}\end{cases}  \\\midrule 2r+1  
&\omega(r,s)  
&\begin{cases}\omega(r,s)\text{\ if\ }r>0\\-1\text{\ otherwise}\end{cases}\end{array}$}

     \item \fbox{$P_{5}  :(a,b)(c,d)=(ac-db^*,da^*+bc)$}\\  
\fbox{$\begin{array}{lll}\omega_{5}&2s&2s+1\\\toprule 2r&\omega(r,s)&\begin{cases}-\omega(s,r)\text{\ if\
}r>0\\1\text{\ otherwise}\end{cases}  \\\midrule 2r+1  
&\omega(r,s)  
&\begin{cases}\omega(s,r)\text{\ if\ }r>0\\-1\text{\ otherwise}\end{cases}\end{array}$}
     \item \fbox{$P_{6}  :(a,b)(c,d)=(ac-b^*d,a^*d+cb)$}\\  
\fbox{$\begin{array}{lll}\omega_{6}&2s&2s+1\\\toprule 2r&\omega(r,s)&\begin{cases}-\omega(r,s)\text{\ if\
}r>0\\1\text{\ otherwise}\end{cases}  \\\midrule 2r+1  
&\omega(s,r)  
&\begin{cases}\omega(r,s)\text{\ if\ }r>0\\-1\text{\ otherwise}\end{cases}\end{array}$}

\newpage

     \item \fbox{$P_{7}  :(a,b)(c,d)=(ac-db^*,a^*d+cb)$}\\  
\fbox{$\begin{array}{lll}\omega_{7}&2s&2s+1\\\toprule 2r&\omega(r,s)&\begin{cases}-\omega(r,s)\text{\ if\
}r>0\\1\text{\ otherwise}\end{cases}  \\\midrule 2r+1  
&\omega(s,r) 
&\begin{cases}\omega(s,r)\text{\ if\ }r>0\\-1\text{\ otherwise}\end{cases}\end{array}$}
     \item \fbox{$P_{8}   :(a,b)(c,d)=(ca-bd^*,da^*+bc)$}\\  
\fbox{$\begin{array}{lll}\omega_{8}&2s&2s+1\\\toprule 2r&\omega(s,r)&\begin{cases}-\omega(s,r)\text{\ if\
}r>0\\1\text{\ otherwise}\end{cases}  \\\midrule 2r+1  
&\omega(r,s)  
&\begin{cases}\omega(r,s)\text{\ if\ }s>0\\-1\text{\ otherwise}\end{cases}\end{array}$}
     \item \fbox{$P_{9}   :(a,b)(c,d)=(ca-d^*b,da^*+bc)$}\\  
\fbox{$\begin{array}{lll}\omega_{9}&2s&2s+1\\\toprule 2r&\omega(s,r)&\begin{cases}-\omega(s,r)\text{\ if\
}r>0\\1\text{\ otherwise}\end{cases}  \\\midrule 2r+1  
&\omega(r,s)  
&\begin{cases}\omega(s,r)\text{\ if\ }s>0\\-1\text{\ otherwise}\end{cases}\end{array}$}
     \item \fbox{$P_{10}   :(a,b)(c,d)=(ca-bd^*,a^*d+cb)$}\\  
\fbox{$\begin{array}{lll}\omega_{10}&2s&2s+1\\\toprule 2r&\omega(s,r)&\begin{cases}-\omega(r,s)\text{\ if\
}r>0\\1\text{\ otherwise}\end{cases}  \\\midrule 2r+1  
&\omega(s,r)  
&\begin{cases}\omega(r,s)\text{\ if\ }s>0\\-1\text{\ otherwise}\end{cases}\end{array}$}

\newpage

     \item \fbox{$P_{11}   :(a,b)(c,d)=(ca-d^*b,a^*d+cb)$}\\  
\fbox{$\begin{array}{lll}\omega_{11}&2s&2s+1\\\toprule 2r&\omega(s,r)&\begin{cases}-\omega(r,s)\text{\ if\
}r>0\\1\text{\ otherwise}\end{cases}  \\\midrule 2r+1  
&\omega(s,r)  
&\begin{cases}\omega(s,r)\text{\ if\ }s>0\\-1\text{\ otherwise}\end{cases}\end{array}$}
     \item \fbox{$P_{12}  :(a,b)(c,d)=(ac-bd^*,da^*+bc)$}\\  
\fbox{$\begin{array}{lll}\omega_{12}&2s&2s+1\\\toprule 2r&\omega(r,s)&\begin{cases}-\omega(s,r)\text{\ if\
}r>0\\1\text{\ otherwise}\end{cases}  \\\midrule 2r+1  
&\omega(r,s)  
&\begin{cases}\omega(r,s)\text{\ if\ }s>0\\-1\text{\ otherwise}\end{cases}\end{array}$}
    \item \fbox{$P_{13}  :(a,b)(c,d)=(ac-d^*b,da^*+bc)$}\\  
\fbox{$\begin{array}{lll}\omega_{13}&2s&2s+1\\\toprule 2r&\omega(r,s)&\begin{cases}-\omega(s,r)\text{\ if\
}r>0\\1\text{\ otherwise}\end{cases}  \\\midrule 2r+1  
&\omega(r,s)  
&\begin{cases}\omega(s,r)\text{\ if\ }s>0\\-1\text{\ otherwise}\end{cases}\end{array}$}
     \item \fbox{$P_{14}  :(a,b)(c,d)=(ac-bd^*,a^*d+cb)$}\\  
\fbox{$\begin{array}{lll}\omega_{14}&2s&2s+1\\\toprule 2r&\omega(r,s)&\begin{cases}-\omega(r,s)\text{\ if\
}r>0\\1\text{\ otherwise}\end{cases}  \\\midrule 2r+1  
&\omega(s,r)  
&\begin{cases}\omega(r,s)\text{\ if\ }s>0\\-1\text{\ otherwise}\end{cases}\end{array}$}

\newpage
 
     \item \fbox{$P_{15}  :(a,b)(c,d)=(ac-d^*b,a^*d+cb)$}\\  
\fbox{$\begin{array}{lll}\omega_{15}&2s&2s+1\\\toprule 2r&\omega(r,s)&\begin{cases}-\omega(r,s)\text{\ if\
}r>0\\1\text{\ otherwise}\end{cases}  \\\midrule 2r+1 
&\omega(s,r) 
&\begin{cases}\omega(s,r)\text{\ if\ }s>0\\-1\text{\ otherwise}\end{cases}\end{array}$}
     \item \fbox{$P_{16}  :(a,b)(c,d)=(ca-b^*d,ad+c^*b)$}\\  
\fbox{$\begin{array}{lll}\omega_{16}&2s&2s+1\\\toprule 2r&\omega(s,r)&\omega(r,s)  \\\midrule 2r+1  
&\begin{cases}-\omega(s,r)\text{\ if\ }s>0\\ 1\text{\ otherwise}\end{cases}  
&\begin{cases}\omega(r,s)\text{\ if\ }r>0\\-1\text{\ otherwise}\end{cases}\end{array}$}
     \item \fbox{$P_{17}   :(a,b)(c,d)=(ca-db^*,ad+c^*b)$}\\  
\fbox{$\begin{array}{lll}\omega_{17}&2s&2s+1\\\toprule 2r&\omega(s,r)&\omega(r,s)  \\\midrule 2r+1  
&\begin{cases}-\omega(s,r)\text{\ if\ }s>0\\ 1\text{\ otherwise}\end{cases}  
&\begin{cases}\omega(s,r)\text{\ if\ }r>0\\-1\text{\ otherwise}\end{cases}\end{array}$}
     \item \fbox{$P_{18}  :(a,b)(c,d)=(ca-b^*d,da+bc^*)$}\\  
\fbox{$\begin{array}{lll}\omega_{18}&2s&2s+1\\\toprule 2r&\omega(s,r)&\omega(s,r)\\\midrule 2r+1  
&\begin{cases}-\omega(r,s)\text{\ if\ }s>0\\ 1\text{\ otherwise}\end{cases}  
&\begin{cases}\omega(r,s)\text{\ if\ }r>0\\-1\text{\ otherwise}\end{cases}\end{array}$}
     \item \fbox{$P_{19}   :(a,b)(c,d)=(ca-db^*,da+bc^*)$}\\  
\fbox{$\begin{array}{lll}\omega_{19}&2s&2s+1\\\toprule 2r&\omega(s,r)&\omega(s,r)\\\midrule 2r+1  
&\begin{cases}-\omega(r,s)\text{\ if\ }s>0\\ 1\text{\ otherwise}\end{cases}  
&\begin{cases}\omega(s,r)\text{\ if\ }r>0\\-1\text{\ otherwise}\end{cases}\end{array}$}

\newpage

     \item \fbox{$P_{20}  :(a,b)(c,d)=(ac-b^*d,ad+c^*b)$}\\  
\fbox{$\begin{array}{lll}\omega_{20}&2s&2s+1\\\toprule 2r&\omega(r,s)&\omega(r,s)  \\\midrule 2r+1  
&\begin{cases}-\omega(s,r)\text{\ if\ }s>0\\ 1\text{\ otherwise}\end{cases}  
&\begin{cases}\omega(r,s)\text{\ if\ }r>0\\-1\text{\ otherwise}\end{cases}\end{array}$}
     \item \fbox{$P_{21}  :(a,b)(c,d)=(ac-db^*,ad+c^*b)$}\\  
\fbox{$\begin{array}{lll}\omega_{21}&2s&2s+1\\\toprule 2r&\omega(r,s)&\omega(r,s)  \\\midrule 2r+1  
&\begin{cases}-\omega(s,r)\text{\ if\ }s>0\\ 1\text{\ otherwise}\end{cases}  
&\begin{cases}\omega(s,r)\text{\ if\ }r>0\\-1\text{\ otherwise}\end{cases}\end{array}$}
     \item \fbox{$P_{22}  :(a,b)(c,d)=(ac-b^*d,da+bc^*)$}\\  
\fbox{$\begin{array}{lll}\omega_{22}&2s&2s+1\\\toprule 2r&\omega(r,s)&\omega(s,r)\\\midrule 2r+1  
&\begin{cases}-\omega(r,s)\text{\ if\ }s>0\\ 1\text{\ otherwise}\end{cases}  
&\begin{cases}\omega(r,s)\text{\ if\ }r>0\\-1\text{\ otherwise}\end{cases}\end{array}$}
     \item \fbox{$P_{23}  :(a,b)(c,d)=(ac-db^*,da+bc^*)$}\\  
\fbox{$\begin{array}{lll}\omega_{23}&2s&2s+1\\\toprule 2r&\omega(r,s)&\omega(s,r)\\\midrule 2r+1  
&\begin{cases}-\omega(r,s)\text{\ if\ }s>0\\ 1\text{\ otherwise}\end{cases}  
&\begin{cases}\omega(s,r)\text{\ if\ }r>0\\-1\text{\ otherwise}\end{cases}\end{array}$}
     \item \fbox{$P_{24}   :(a,b)(c,d)=(ca-bd^*,ad+c^*b)$}\\  
\fbox{$\begin{array}{lll}\omega_{24}&2s&2s+1\\\toprule 2r&\omega(s,r)&\omega(r,s)  \\\midrule 2r+1  
&\begin{cases}-\omega(s,r)\text{\ if\ }s>0\\ 1\text{\ otherwise}\end{cases}  
&\begin{cases}\omega(r,s)\text{\ if\ }s>0\\-1\text{\ otherwise}\end{cases}\end{array}$}

\newpage

     \item \fbox{$P_{25}   :(a,b)(c,d)=(ca-d^*b,ad+c^*b)$}\\  
\fbox{$\begin{array}{lll}\omega_{25}&2s&2s+1\\\toprule 2r&\omega(s,r)&\omega(r,s)  \\\midrule 2r+1  
&\begin{cases}-\omega(s,r)\text{\ if\ }s>0\\ 1\text{\ otherwise}\end{cases}  
&\begin{cases}\omega(s,r)\text{\ if\ }s>0\\-1\text{\ otherwise}\end{cases}\end{array}$}
     \item \fbox{$P_{26}   :(a,b)(c,d)=(ca-bd^*,da+bc^*)$}\\  
\fbox{$\begin{array}{lll}\omega_{26}&2s&2s+1\\\toprule 2r&\omega(s,r)&\omega(s,r)\\\midrule 2r+1  
&\begin{cases}-\omega(r,s)\text{\ if\ }s>0\\ 1\text{\ otherwise}\end{cases}  
&\begin{cases}\omega(r,s)\text{\ if\ }s>0\\-1\text{\ otherwise}\end{cases}\end{array}$}
     \item \fbox{$P_{27}   :(a,b)(c,d)=(ca-d^*b,da+bc^*)$}\\  
\fbox{$\begin{array}{lll}\omega_{27}&2s&2s+1\\\toprule 2r&\omega(s,r)&\omega(s,r)\\\midrule 2r+1  
&\begin{cases}-\omega(r,s)\text{\ if\ }s>0\\ 1\text{\ otherwise}\end{cases}
&\begin{cases}\omega(s,r)\text{\ if\ }s>0\\-1\text{\ otherwise}\end{cases}\end{array}$}
     \item \fbox{$P_{28}  :(a,b)(c,d)=(ac-bd^*,ad+c^*b)$}\\  
\fbox{$\begin{array}{lll}\omega_{28}&2s&2s+1\\\toprule 2r&\omega(r,s)&\omega(r,s)  \\\midrule 2r+1  
&\begin{cases}-\omega(s,r)\text{\ if\ }s>0\\ 1\text{\ otherwise}\end{cases}  
&\begin{cases}\omega(r,s)\text{\ if\ }s>0\\-1\text{\ otherwise}\end{cases}\end{array}$}
     \item \fbox{$P_{29}  :(a,b)(c,d)=(ac-d^*b,ad+c^*b)$}\\  
\fbox{$\begin{array}{lll}\omega_{29}&2s&2s+1\\\toprule 2r&\omega(r,s)&\omega(r,s)  \\\midrule 2r+1  
&\begin{cases}-\omega(s,r)\text{\ if\ }s>0\\ 1\text{\ otherwise}\end{cases}  
&\begin{cases}\omega(s,r)\text{\ if\ }s>0\\-1\text{\ otherwise}\end{cases}\end{array}$}

\newpage

     \item \fbox{$P_{30}  :(a,b)(c,d)=(ac-bd^*,da+bc^*)$}\\  
\fbox{$\begin{array}{lll}\omega_{30}&2s&2s+1\\\toprule 2r&\omega(r,s)&\omega(s,r)\\\midrule 2r+1  
&\begin{cases}-\omega(r,s)\text{\ if\ }s>0\\ 1\text{\ otherwise}\end{cases}  
&\begin{cases}\omega(r,s)\text{\ if\ }s>0\\-1\text{\ otherwise}\end{cases}\end{array}$}
     \item \fbox{$P_{31}  :(a,b)(c,d)=(ac-d^*b,da+bc^*)$}\\  
\fbox{$\begin{array}{lll}\omega_{31}&2s&2s+1\\\toprule 2r&\omega(r,s)&\omega(s,r)\\\midrule 2r+1  
&\begin{cases}-\omega(r,s)\text{\ if\ }s>0\\ 1\text{\ otherwise}\end{cases}  
&\begin{cases}\omega(s,r)\text{\ if\ }s>0\\-1\text{\ otherwise}\end{cases}\end{array}$}
\end{enumerate}

\section{The twist blocks}

For each of the 32 twists $\omega_k$ on the group $\{0,1,2,3,\cdots\}$ there are five fundamental
$2\times2$ matrices
$E_{krs}=\left[\begin{array}{ll}\omega_k(2r,2s)&\omega_k(2r,2s+1)\\\omega_k(2r+1,2s)&\omega_k(2r+1,2s+1)\end{array}
\right ] $ which make up the twist table. 
\begin{itemize}
  \item For $r=s=0$, $E_{k00}=\omega_k(0,0)C$, where $C$ is the block in the upper left corner of the twist table.
  \item For $r>s=0$, $E_{kr0}=\omega_k(r,0)L$, where $L$ is the block occupying every position in the left edge of the
twist table below the corner block.
  \item For $s>r=0$, $E_{k0s}=\omega_k(0,s)T$, where $T$ is the block occupying every position in the top edge of the
twist table to the right of the corner block.
  \item For $r=s\ne0$, $E_{krr}=\omega_k(r,r)D$, where $D$ is the block whose negative occupies every position on the
diagonal of the twist table with the exception of the corner block.
  \item For $0\ne r\ne s\ne 0$, $E_{krs}=\omega_k(r,s)N$, where $N$ is the interior block of the twist table. Either
$N$ or its negative (depending upon the value of $\omega(r,s)$) occupies the twist table in block position $(r,s)$
(counting the first block row as block row zero and the first block column as block column 0) when $0\ne r\ne s\ne 0$.
\end{itemize}

The values of the five fundamental blocks $\mathbf{C}, \mathbf{L}, \mathbf{T}, \mathbf{D}$ and $\mathbf{N}$
independent of the value of $\omega(r,s)$ are given for each of the 32 twists in Table \ref{Tfirst} and Table
\ref{Tsecond}.

\begin{table}
\begin{tabular}{cccccc}
 Twist  & $\mathbf{C}$  & $\mathbf{L}$ & $\mathbf{T}$   & $\mathbf{D}$   & $\mathbf{N}$ \\\toprule
        & $r=s=0$ & $r>s=0$   & $s>r=0$    & $r=s\ne0$  & $ 0\ne r\ne s\ne 0$ \\\midrule
 $\omega_{0}$ & $\left[ \begin{array}{rr}  1 &  1\\  1 & -1 \end{array} \right]$
              & $\left[ \begin{array}{rr}  1 & -1\\  1 &  1 \end{array} \right]$
              & $\left[ \begin{array}{rr}  1 &  1\\  1 & -1 \end{array} \right]$
              & $\left[ \begin{array}{rr}  1 & -1\\  1 &  1 \end{array} \right]$
              & $\left[ \begin{array}{rr} -1 &  1\\  1 &  1 \end{array} \right]$\\\midrule

 $\omega_{1}$ & $\left[ \begin{array}{rr}  1 &  1\\  1 & -1 \end{array} \right]$
              & $\left[ \begin{array}{rr}  1 & -1\\  1 &  1 \end{array} \right]$
              & $\left[ \begin{array}{rr}  1 &  1\\  1 & -1 \end{array} \right]$
              & $\left[ \begin{array}{rr}  1 & -1\\  1 &  1 \end{array} \right]$
              & $\left[ \begin{array}{rr} -1 &  1\\  1 & -1 \end{array} \right]$\\\midrule

 $\omega_{2}$ & $\left[ \begin{array}{rr}  1 &  1\\  1 & -1 \end{array} \right]$
              & $\left[ \begin{array}{rr}  1 & -1\\  1 &  1 \end{array} \right]$
              & $\left[ \begin{array}{rr}  1 &  1\\  1 & -1 \end{array} \right]$
              & $\left[ \begin{array}{rr}  1 & -1\\  1 &  1 \end{array} \right]$
              & $\left[ \begin{array}{rr} -1 & -1\\ -1 &  1 \end{array} \right]$\\\bottomrule

 $\omega_{3}$ & $\left[ \begin{array}{rr}  1 &  1\\  1 & -1 \end{array} \right]$
              & $\left[ \begin{array}{rr}  1 & -1\\  1 &  1 \end{array} \right]$
              & $\left[ \begin{array}{rr}  1 &  1\\  1 & -1 \end{array} \right]$
              & $\left[ \begin{array}{rr}  1 & -1\\  1 &  1 \end{array} \right]$
              & $\left[ \begin{array}{rr} -1 & -1\\ -1 & -1 \end{array} \right]$\\\midrule

 $\omega_{4}$ & $\left[ \begin{array}{rr}  1 &  1\\  1 & -1 \end{array} \right]$
              & $\left[ \begin{array}{rr}  1 & -1\\  1 &  1 \end{array} \right]$
              & $\left[ \begin{array}{rr}  1 &  1\\  1 & -1 \end{array} \right]$
              & $\left[ \begin{array}{rr}  1 & -1\\  1 &  1 \end{array} \right]$
              & $\left[ \begin{array}{rr}  1 &  1\\  1 &  1 \end{array} \right]$\\\midrule

 $\omega_{5}$ & $\left[ \begin{array}{rr}  1 &  1\\  1 & -1 \end{array} \right]$
              & $\left[ \begin{array}{rr}  1 & -1\\  1 &  1 \end{array} \right]$
              & $\left[ \begin{array}{rr}  1 &  1\\  1 & -1 \end{array} \right]$
              & $\left[ \begin{array}{rr}  1 & -1\\  1 &  1 \end{array} \right]$
              & $\left[ \begin{array}{rr}  1 &  1\\  1 & -1 \end{array} \right]$\\\midrule

 $\omega_{6}$ & $\left[ \begin{array}{rr}  1 &  1\\  1 & -1 \end{array} \right]$
              & $\left[ \begin{array}{rr}  1 & -1\\  1 &  1 \end{array} \right]$
              & $\left[ \begin{array}{rr}  1 &  1\\  1 & -1 \end{array} \right]$
              & $\left[ \begin{array}{rr}  1 & -1\\  1 &  1 \end{array} \right]$
              & $\left[ \begin{array}{rr}  1 & -1\\ -1 &  1 \end{array} \right]$\\\bottomrule

 $\omega_{7}$ & $\left[ \begin{array}{rr}  1 &  1\\  1 & -1 \end{array} \right]$
              & $\left[ \begin{array}{rr}  1 & -1\\  1 &  1 \end{array} \right]$
              & $\left[ \begin{array}{rr}  1 &  1\\  1 & -1 \end{array} \right]$
              & $\left[ \begin{array}{rr}  1 & -1\\  1 &  1 \end{array} \right]$
              & $\left[ \begin{array}{rr}  1 & -1\\ -1 & -1 \end{array} \right]$\\\midrule

 $\omega_{8}$ & $\left[ \begin{array}{rr}  1 &  1\\  1 & -1 \end{array} \right]$
              & $\left[ \begin{array}{rr}  1 & -1\\  1 & -1 \end{array} \right]$
              & $\left[ \begin{array}{rr}  1 &  1\\  1 &  1 \end{array} \right]$
              & $\left[ \begin{array}{rr}  1 & -1\\  1 &  1 \end{array} \right]$
              & $\left[ \begin{array}{rr} -1 &  1\\  1 &  1 \end{array} \right]$\\\midrule

 $\omega_{9}$ & $\left[ \begin{array}{rr}  1 &  1\\  1 & -1 \end{array} \right]$
              & $\left[ \begin{array}{rr}  1 & -1\\  1 & -1 \end{array} \right]$
              & $\left[ \begin{array}{rr}  1 &  1\\  1 &  1 \end{array} \right]$
              & $\left[ \begin{array}{rr}  1 & -1\\  1 &  1 \end{array} \right]$
              & $\left[ \begin{array}{rr} -1 &  1\\  1 & -1 \end{array} \right]$\\\midrule

$\omega_{10}$ & $\left[ \begin{array}{rr}  1 &  1\\  1 & -1 \end{array} \right]$
              & $\left[ \begin{array}{rr}  1 & -1\\  1 & -1 \end{array} \right]$
              & $\left[ \begin{array}{rr}  1 &  1\\  1 &  1 \end{array} \right]$
              & $\left[ \begin{array}{rr}  1 & -1\\  1 &  1 \end{array} \right]$
              & $\left[ \begin{array}{rr} -1 & -1\\ -1 &  1 \end{array} \right]$\\\midrule

$\omega_{11}$ & $\left[ \begin{array}{rr}  1 &  1\\  1 & -1 \end{array} \right]$
              & $\left[ \begin{array}{rr}  1 & -1\\  1 & -1 \end{array} \right]$
              & $\left[ \begin{array}{rr}  1 &  1\\  1 &  1 \end{array} \right]$
              & $\left[ \begin{array}{rr}  1 & -1\\  1 &  1 \end{array} \right]$
              & $\left[ \begin{array}{rr} -1 & -1\\ -1 & -1 \end{array} \right]$\\\midrule

$\omega_{12}$ & $\left[ \begin{array}{rr}  1 &  1\\  1 & -1 \end{array} \right]$
              & $\left[ \begin{array}{rr}  1 & -1\\  1 & -1 \end{array} \right]$
              & $\left[ \begin{array}{rr}  1 &  1\\  1 &  1 \end{array} \right]$
              & $\left[ \begin{array}{rr}  1 & -1\\  1 &  1 \end{array} \right]$
              & $\left[ \begin{array}{rr}  1 &  1\\  1 &  1 \end{array} \right]$\\\midrule

$\omega_{13}$ & $\left[ \begin{array}{rr}  1 &  1\\  1 & -1 \end{array} \right]$
              & $\left[ \begin{array}{rr}  1 & -1\\  1 & -1 \end{array} \right]$
              & $\left[ \begin{array}{rr}  1 &  1\\  1 &  1 \end{array} \right]$
              & $\left[ \begin{array}{rr}  1 & -1\\  1 &  1 \end{array} \right]$
              & $\left[ \begin{array}{rr}  1 &  1\\  1 & -1 \end{array} \right]$\\\midrule

$\omega_{14}$ & $\left[ \begin{array}{rr}  1 &  1\\  1 & -1 \end{array} \right]$
              & $\left[ \begin{array}{rr}  1 & -1\\  1 & -1 \end{array} \right]$
              & $\left[ \begin{array}{rr}  1 &  1\\  1 &  1 \end{array} \right]$
              & $\left[ \begin{array}{rr}  1 & -1\\  1 &  1 \end{array} \right]$
              & $\left[ \begin{array}{rr}  1 & -1\\ -1 &  1 \end{array} \right]$\\\midrule

$\omega_{15}$ & $\left[ \begin{array}{rr}  1 &  1\\  1 & -1 \end{array} \right]$
              & $\left[ \begin{array}{rr}  1 & -1\\  1 & -1 \end{array} \right]$
              & $\left[ \begin{array}{rr}  1 &  1\\  1 &  1 \end{array} \right]$
              & $\left[ \begin{array}{rr}  1 & -1\\  1 &  1 \end{array} \right]$
              & $\left[ \begin{array}{rr}  1 & -1\\ -1 & -1 \end{array} \right]$\\\midrule

\end{tabular}\\[12pt]
  \caption{Twist blocks for $\omega_{0}$ through $\omega_{15}$}
\label{Tfirst}
\end{table}

\begin{table}
\begin{tabular}{cccccc}
 Twist  & $\mathbf{C}$  & $\mathbf{L}$ & $\mathbf{T}$   & $\mathbf{D}$   & $\mathbf{N}$ \\\toprule
        & $r=s=0$ & $r>s=0$   & $s>r=0$    & $r=s\ne0$  & $ 0\ne r\ne s\ne 0$ \\\midrule
$\omega_{16}$ & $\left[ \begin{array}{rr}  1 &  1\\  1 & -1 \end{array} \right]$
              & $\left[ \begin{array}{rr}  1 &  1\\  1 &  1 \end{array} \right]$
              & $\left[ \begin{array}{rr}  1 &  1\\ -1 & -1 \end{array} \right]$
              & $\left[ \begin{array}{rr}  1 &  1\\ -1 &  1 \end{array} \right]$
              & $\left[ \begin{array}{rr} -1 &  1\\  1 &  1 \end{array} \right]$\\\midrule

$\omega_{17}$ & $\left[ \begin{array}{rr}  1 &  1\\  1 & -1 \end{array} \right]$
              & $\left[ \begin{array}{rr}  1 &  1\\  1 &  1 \end{array} \right]$
              & $\left[ \begin{array}{rr}  1 &  1\\ -1 & -1 \end{array} \right]$
              & $\left[ \begin{array}{rr}  1 &  1\\ -1 &  1 \end{array} \right]$
              & $\left[ \begin{array}{rr} -1 &  1\\  1 & -1 \end{array} \right]$\\\midrule

$\omega_{18}$ & $\left[ \begin{array}{rr}  1 &  1\\  1 & -1 \end{array} \right]$
              & $\left[ \begin{array}{rr}  1 &  1\\  1 &  1 \end{array} \right]$
              & $\left[ \begin{array}{rr}  1 &  1\\ -1 & -1 \end{array} \right]$
              & $\left[ \begin{array}{rr}  1 &  1\\ -1 &  1 \end{array} \right]$
              & $\left[ \begin{array}{rr} -1 & -1\\ -1 &  1 \end{array} \right]$\\\midrule

$\omega_{19}$ & $\left[ \begin{array}{rr}  1 &  1\\  1 & -1 \end{array} \right]$
              & $\left[ \begin{array}{rr}  1 &  1\\  1 &  1 \end{array} \right]$
              & $\left[ \begin{array}{rr}  1 &  1\\ -1 & -1 \end{array} \right]$
              & $\left[ \begin{array}{rr}  1 &  1\\ -1 &  1 \end{array} \right]$
              & $\left[ \begin{array}{rr} -1 & -1\\ -1 & -1 \end{array} \right]$\\\midrule

$\omega_{20}$ & $\left[ \begin{array}{rr}  1 &  1\\  1 & -1 \end{array} \right]$
              & $\left[ \begin{array}{rr}  1 &  1\\  1 &  1 \end{array} \right]$
              & $\left[ \begin{array}{rr}  1 &  1\\ -1 & -1 \end{array} \right]$
              & $\left[ \begin{array}{rr}  1 &  1\\ -1 &  1 \end{array} \right]$
              & $\left[ \begin{array}{rr}  1 &  1\\  1 &  1 \end{array} \right]$\\\midrule

$\omega_{21}$ & $\left[ \begin{array}{rr}  1 &  1\\  1 & -1 \end{array} \right]$
              & $\left[ \begin{array}{rr}  1 &  1\\  1 &  1 \end{array} \right]$
              & $\left[ \begin{array}{rr}  1 &  1\\ -1 & -1 \end{array} \right]$
              & $\left[ \begin{array}{rr}  1 &  1\\ -1 &  1 \end{array} \right]$
              & $\left[ \begin{array}{rr}  1 &  1\\  1 & -1 \end{array} \right]$\\\midrule

$\omega_{22}$ & $\left[ \begin{array}{rr}  1 &  1\\  1 & -1 \end{array} \right]$
              & $\left[ \begin{array}{rr}  1 &  1\\  1 &  1 \end{array} \right]$
              & $\left[ \begin{array}{rr}  1 &  1\\ -1 & -1 \end{array} \right]$
              & $\left[ \begin{array}{rr}  1 &  1\\ -1 &  1 \end{array} \right]$
              & $\left[ \begin{array}{rr}  1 & -1\\ -1 &  1 \end{array} \right]$\\\midrule

$\omega_{23}$ & $\left[ \begin{array}{rr}  1 &  1\\  1 & -1 \end{array} \right]$
              & $\left[ \begin{array}{rr}  1 &  1\\  1 &  1 \end{array} \right]$
              & $\left[ \begin{array}{rr}  1 &  1\\ -1 & -1 \end{array} \right]$
              & $\left[ \begin{array}{rr}  1 &  1\\ -1 &  1 \end{array} \right]$
              & $\left[ \begin{array}{rr}  1 & -1\\ -1 & -1 \end{array} \right]$\\\midrule

$\omega_{24}$ & $\left[ \begin{array}{rr}  1 &  1\\  1 & -1 \end{array} \right]$
              & $\left[ \begin{array}{rr}  1 &  1\\  1 & -1 \end{array} \right]$
              & $\left[ \begin{array}{rr}  1 &  1\\ -1 &  1 \end{array} \right]$
              & $\left[ \begin{array}{rr}  1 &  1\\ -1 &  1 \end{array} \right]$
              & $\left[ \begin{array}{rr} -1 &  1\\  1 &  1 \end{array} \right]$\\\midrule

$\omega_{25}$ & $\left[ \begin{array}{rr}  1 &  1\\  1 & -1 \end{array} \right]$
              & $\left[ \begin{array}{rr}  1 &  1\\  1 & -1 \end{array} \right]$
              & $\left[ \begin{array}{rr}  1 &  1\\ -1 &  1 \end{array} \right]$
              & $\left[ \begin{array}{rr}  1 &  1\\ -1 &  1 \end{array} \right]$
              & $\left[ \begin{array}{rr} -1 &  1\\  1 & -1 \end{array} \right]$\\\midrule

$\omega_{26}$ & $\left[ \begin{array}{rr}  1 &  1\\  1 & -1 \end{array} \right]$
              & $\left[ \begin{array}{rr}  1 &  1\\  1 & -1 \end{array} \right]$
              & $\left[ \begin{array}{rr}  1 &  1\\ -1 &  1 \end{array} \right]$
              & $\left[ \begin{array}{rr}  1 &  1\\ -1 &  1 \end{array} \right]$
              & $\left[ \begin{array}{rr} -1 & -1\\ -1 &  1 \end{array} \right]$\\\midrule

$\omega_{27}$ & $\left[ \begin{array}{rr}  1 &  1\\  1 & -1 \end{array} \right]$
              & $\left[ \begin{array}{rr}  1 &  1\\  1 & -1 \end{array} \right]$
              & $\left[ \begin{array}{rr}  1 &  1\\ -1 &  1 \end{array} \right]$
              & $\left[ \begin{array}{rr}  1 &  1\\ -1 &  1 \end{array} \right]$
              & $\left[ \begin{array}{rr} -1 & -1\\ -1 & -1 \end{array} \right]$\\\midrule
$\omega_{28}$ & $\left[ \begin{array}{rr}  1 &  1\\  1 & -1 \end{array} \right]$
              & $\left[ \begin{array}{rr}  1 &  1\\  1 & -1 \end{array} \right]$
              & $\left[ \begin{array}{rr}  1 &  1\\ -1 &  1 \end{array} \right]$
              & $\left[ \begin{array}{rr}  1 &  1\\ -1 &  1 \end{array} \right]$
              & $\left[ \begin{array}{rr}  1 &  1\\  1 &  1 \end{array} \right]$\\\midrule

$\omega_{29}$ & $\left[ \begin{array}{rr}  1 &  1\\  1 & -1 \end{array} \right]$
              & $\left[ \begin{array}{rr}  1 &  1\\  1 & -1 \end{array} \right]$
              & $\left[ \begin{array}{rr}  1 &  1\\ -1 &  1 \end{array} \right]$
              & $\left[ \begin{array}{rr}  1 &  1\\ -1 &  1 \end{array} \right]$
              & $\left[ \begin{array}{rr}  1 &  1\\  1 & -1 \end{array} \right]$\\\midrule

$\omega_{30}$ & $\left[ \begin{array}{rr}  1 &  1\\  1 & -1 \end{array} \right]$
              & $\left[ \begin{array}{rr}  1 &  1\\  1 & -1 \end{array} \right]$
              & $\left[ \begin{array}{rr}  1 &  1\\ -1 &  1 \end{array} \right]$
              & $\left[ \begin{array}{rr}  1 &  1\\ -1 &  1 \end{array} \right]$
              & $\left[ \begin{array}{rr}  1 & -1\\ -1 &  1 \end{array} \right]$\\\midrule

$\omega_{31}$ & $\left[ \begin{array}{rr}  1 &  1\\  1 & -1 \end{array} \right]$
              & $\left[ \begin{array}{rr}  1 &  1\\  1 & -1 \end{array} \right]$
              & $\left[ \begin{array}{rr}  1 &  1\\ -1 &  1 \end{array} \right]$
              & $\left[ \begin{array}{rr}  1 &  1\\ -1 &  1 \end{array} \right]$
              & $\left[ \begin{array}{rr}  1 & -1\\ -1 & -1 \end{array} \right]$\\\midrule

\end{tabular}\\[12pt]

  \caption{Twist blocks for $\omega_{16}$ through $\omega_{31}$}
  \label{Tsecond}
\end{table}

\newpage

\section{The Twist Tree}

Looking at tables \ref{Tfirst} and \ref{Tsecond} it is seen that the five basic twist blocks are:

\begin{equation}\label{Cblock}
  \mathbf{C}=\left[\begin{array}{rr} 1 & 1 \\ 1 & -1\end{array}\right]
\end{equation}

\begin{equation}\label{Lblock}
  \mathbf{L}=\left[\begin{array}{rr} 1 & \alpha^\prime \\ 1 & \alpha\end{array}\right]
\end{equation}

\begin{equation}\label{Tblock}
  \mathbf{T}=\left[\begin{array}{rr} 1 & 1 \\ \alpha & \alpha^\prime\end{array}\right]
\end{equation}

\begin{equation}\label{Dblock}
  \mathbf{D}=\left[\begin{array}{rr} 1 & \alpha^\prime \\ \alpha & 1\end{array}\right]
\end{equation}

\begin{equation}\label{Nblock}
  \mathbf{N}=\left[\begin{array}{rr} \beta & \gamma \\ \gamma & \gamma^\prime\end{array}\right]
\end{equation}

The five twist constants $\alpha^\prime,\alpha,\beta,\gamma,\gamma^\prime$ are always equal to either $1$ or $-1$. 

Square matrices having $2^N$ rows and columns can be partitioned into $2\times2$ partitioned block matrices as follows:

\begin{equation}\label{E:pmat}
M = \begin{pmat}[{|}]
              M_{0,0} & M_{0,1} \cr\-
              M_{1,0} & M_{1,1}\cr
        \end{pmat}
\end{equation}

Each submatrix, provided it is not a $1\times 1$ matrix, can be further subdivided. For example\\

\begin{equation}
M_{1,0} = \begin{pmat}[{|}]
              M_{10,00} & M_{10,01} \cr\-
              M_{11,00} & M_{11,01}\cr
        \end{pmat}
\end{equation}

Then if the rows and columns of the matrices are numbered beginning with 0 rather than 1 and when, ultimately $M_{p,q}$ is a $1\times1$ matrix, it
will be the case that

\begin{equation}
  M_{i,j}=m_{ij}
\end{equation}

Equation \ref{E:pmat} can be written in `tree' form as follows:

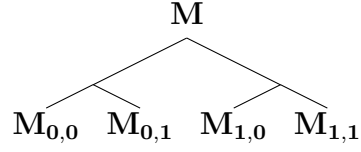
\Tree [.$\mathbf{M}$ [. $\mathbf{ M_{0,0}}$ $\mathbf{ M_{0,1}}$ ] [. $\mathbf{M_{1,0} }$ $\mathbf{M_{1,1} }$ ] ]\\

The subscripts can be viewed as navigation instructions: A \emph{zero} is an instruction to move down a left branch and a \emph{one} is an instruction
to move down a right branch of the tree. Each of the block matrices $\mathbf{C},\mathbf{L},\mathbf{T},\mathbf{D}$ and $\mathbf{N}$ can be written in
tree form.

The twist tree for $\omega_{7}$ was derived in \cite{B2009}. Table \ref{fig:cydtreeG} shows the general twist tree
for all 32 variations on the Cayley-Dickson doubling formula. 

\begin{table}[ht]
~\\
\Tree [.$\mathbf{C}$ [. $\mathbf{C}$ $\mathbf{T}$ ] [. $\mathbf{L}$ $\mathbf{-D}$ ] ]
\Tree [.$\mathbf{L}$ [. $\mathbf{L}$ $\alpha^\prime\mathbf{N}$ ] [. $\mathbf{L}$ $\alpha\mathbf{N}$ ] ]\\[12pt]
\Tree [.$\mathbf{T}$ [. $\mathbf{T}$ $\mathbf{T}$ ] [. $-\alpha^\prime\mathbf{N}$ $-\alpha\mathbf{N}$ ] ]
\Tree [.$\mathbf{-D}$ [. $\mathbf{-D}$ $-\alpha^\prime\mathbf{N}$ ] [ $\alpha^\prime\mathbf{N}$ $\mathbf{-D}$ ] ]\\[12pt]
\Tree [.$\mathbf{N}$ [. $\beta\mathbf{N}$ $\gamma\mathbf{N}$ ] [. $\gamma\mathbf{N}$ $\gamma^\prime\mathbf{N}$ ] ]
\Tree [.$-\mathbf{N}$ [. $-\beta\mathbf{N}$ $-\gamma\mathbf{N}$ ] [. $-\gamma\mathbf{N}$ $-\gamma^\prime\mathbf{N}$ ] ] 
\vspace{3mm}
\caption{Generalized `Cayley-Dickson-Like' Twist Tree}
\label{fig:cydtreeG}
\end{table}

The values of $\alpha,\alpha^\prime,\beta,\gamma$ and $\gamma^\prime$ for the 32 different doublings are given in table \vref{Tthird}.

\begin{table}[ht]
\begin{tabular}{l|ccccc}
& $\alpha^\prime$ & $\alpha$ & $\beta$ & $\gamma$ & $\gamma^\prime$\\\toprule
$P_{0}:   (a,b)(c,d)=(ca-b^*d,da^*+bc)$&$-$ &$+$ &$-$ &$+$ & $+$ \\  
$P_{1}:   (a,b)(c,d)=(ca-db^*,da^*+bc)$&$-$ &$+$ &$-$ &$+$ &$-$ \\  
$P_{2}:   (a,b)(c,d)=(ca-b^*d,a^*d+cb)$&$-$ &$+$ &$-$ &$-$&$+$ \\  
$P_{3}:   (a,b)(c,d)=(ca-db^*,a^*d+cb)$&$-$ &$+$ &$-$ &$-$&$-$ \\  
$P_{4}:   (a,b)(c,d)=(ac-b^*d,da^*+bc)$&$-$ &$+$ &$+$ &$+$ &$+$ \\  
$P_{5}:   (a,b)(c,d)=(ac-db^*,da^*+bc)$&$-$ &$+$ &$+$ &$+$ &$-$ \\  
$P_{6}:   (a,b)(c,d)=(ac-b^*d,a^*d+cb)$&$-$ &$+$ &$+$ &$-$ &$+$ \\  
$P_{7}:   (a,b)(c,d)=(ac-db^*,a^*d+cb)$&$-$ &$+$ &$+$ &$-$ &$-$ \\  
$P_{8}:   (a,b)(c,d)=(ca-bd^*,da^*+bc)$&$-$ &$-$&$-$ &$+$ &$+$ \\  
$P_{9}:   (a,b)(c,d)=(ca-d^*b,da^*+bc)$&$-$ &$-$&$-$ &$+$ &$-$ \\  
$P_{10}:  (a,b)(c,d)=(ca-bd^*,a^*d+cb)$&$-$ &$-$&$-$ &$-$&$+$ \\  
$P_{11}:  (a,b)(c,d)=(ca-d^*b,a^*d+cb)$&$-$ &$-$&$-$ &$-$&$-$ \\  
$P_{12}:  (a,b)(c,d)=(ac-bd^*,da^*+bc)$&$-$ &$-$&$+$ &$+$ &$+$ \\  
$P_{13}:  (a,b)(c,d)=(ac-d^*b,da^*+bc)$&$-$ &$-$&$+$ &$+$ &$-$ \\  
$P_{14}:  (a,b)(c,d)=(ac-bd^*,a^*d+cb)$&$-$ &$-$&$+$ &$-$ &$+$ \\  
$P_{15}:  (a,b)(c,d)=(ac-d^*b,a^*d+cb)$&$-$ &$-$&$+$ &$-$ &$-$ \\  
$P_{16}:  (a,b)(c,d)=(ca-b^*d,ad+c^*b)$&$+$ &$+$ &$-$ &$+$ &$+$ \\  
$P_{17}:  (a,b)(c,d)=(ca-db^*,ad+c^*b)$&$+$ &$+$ &$-$ &$+$ &$-$ \\  
$P_{18}:  (a,b)(c,d)=(ca-b^*d,da+bc^*)$&$+$ &$+$ &$-$ &$-$&$+$ \\  
$P_{19}:  (a,b)(c,d)=(ca-db^*,da+bc^*)$&$+$ &$+$ &$-$ &$-$&$-$ \\  
$P_{20}:  (a,b)(c,d)=(ac-b^*d,ad+c^*b)$&$+$ &$+$ &$+$ &$+$ &$+$ \\  
$P_{21}:  (a,b)(c,d)=(ac-db^*,ad+c^*b)$&$+$ &$+$ &$+$ &$+$ &$-$ \\  
$P_{22}:  (a,b)(c,d)=(ac-b^*d,da+bc^*)$&$+$ &$+$ &$+$ &$-$ &$+$ \\  
$P_{23}:  (a,b)(c,d)=(ac-db^*,da+bc^*)$&$+$ &$+$ &$+$ &$-$ &$-$ \\  
$P_{24}:  (a,b)(c,d)=(ca-bd^*,ad+c^*b)$&$+$ &$-$ &$-$&$+$ &$+$ \\  
$P_{25}:  (a,b)(c,d)=(ca-d^*b,ad+c^*b)$&$+$ &$-$ &$-$&$+$ &$-$ \\  
$P_{26}:  (a,b)(c,d)=(ca-bd^*,da+bc^*)$&$+$ &$-$ &$-$&$-$ &$+$ \\  
$P_{27}:  (a,b)(c,d)=(ca-d^*b,da+bc^*)$&$+$ &$-$ &$-$&$-$ &$-$ \\  
$P_{28}:  (a,b)(c,d)=(ac-bd^*,ad+c^*b)$&$+$ &$-$ &$+$ &$+$ &$+$ \\  
$P_{29}:  (a,b)(c,d)=(ac-d^*b,ad+c^*b)$&$+$ &$-$ &$+$ &$+$ &$-$ \\  
$P_{30}:  (a,b)(c,d)=(ac-bd^*,da+bc^*)$&$+$ &$-$ &$+$ &$-$ &$+$ \\  
$P_{31}:  (a,b)(c,d)=(ac-d^*b,da+bc^*)$&$+$ &$-$ &$+$ &$-$ &$-$ \\\bottomrule  
\end{tabular}\\[12pt]
\caption{Constants for the 32 doubling products}
\label{Tthird}
\end{table}

\section{The quaternion properties}

The two quaternion properties, stated in terms of the twist $\omega$ are as follows:\\

If $0\ne p\ne q\ne0$ then

\begin{enumerate}
	\item $\omega(p,q)+\omega(q,p)=0$
	\item $\omega(p,q)=\omega(q,pq)=\omega(pq,p)$
\end{enumerate}

$\omega(q,p)+\omega(p,q)=0$ was previously shown to be true for all 32 products (see equation \vref{Eqn:neg}).

It will be shown that $\omega(p,q)=\omega(q,pq)=\omega(pq,p)$ for only eight of the 32 Cayley-Dickson-like doubling products. These eight are the
only actual Cayley-Dickson products.

\subsection{The basis step of the induction}

Begin by showing that the property holds when either $p$ or $q$ equals 1 or when $p$ and $q$ differ by 1. This forms the basis step of an induction.

For $p=2r\ne 0$ and $q=1$ we have $pq=2r+1,$ so we must compare the values of $\omega(2r,1),\omega(1,2r+1)$ and $\omega(2r+1,2r).$

$\omega(2r,1)=\omega(2r,2(0)+1)$ so we find the value in the $(2r,2s+1)$ portions of the charts with $r>0,s=0.$
We see that the value is $-1$ for $\omega_{0}$ through $\omega_{15}$ and 1 for $\omega_{16}$ through $\omega_{31}.$

$\omega(1,2r+1)=\omega(2(0)+1,2s+1)$ We find these values in the $(2r+1,2s+1)$ portion of the tables with $r=0$ and $s>0.$  We see that the value is $-1$ for $\omega_{0}$ through $\omega_{7}$ and $\omega_{16}$ through $\omega_{23}$ and 1 for $\omega_{8}$ through $\omega_{15}$ and $\omega_{24}$ through $\omega_{31}$.

 We conclude that  the twists $\omega_{8}$ through $\omega_{23}$ fail to have the second quaternion property.
 
 $\omega(2r,2r+1)$ has a value of $-1$ for  $\omega_{0}$ through $\omega_{15}$ and 1 for  $\omega_{16}$ through $\omega_{31}.$ 
 
 Thus $\omega(2r,1)=\omega(1,2r+1)=\omega(2r+1,2r)$ for  $\omega_{0}$ through $\omega_{7}$ and for $\omega_{24}$ through $\omega_{31}$ alone.
 
 Since $\omega(2r,1)=-\omega(1,2r),\omega(1,2r+1)=-\omega(2r+1,1)$ and $\omega(2r+1,2r)=-\omega(2r,2r+1)$ we have the result that $\omega(2r+1,1)=\omega(1,2r)=\omega(2r,2r+1).$
 
 So for $\omega_{0}$ through $\omega_{7}$ and for $\omega_{24}$ through $\omega_{31}$ alone it is the case that $\omega(p,q)=\omega(q,pq)=\omega(pq,p)$ in the cases where
 
 \begin{enumerate}
 	\item $p=2r>0,q=1$
 	\item $p=2r+1,q=1$
 	\item $p=2r>0,q=2r+1$
 	\item $p=2r+1,q=2r+1$
 \end{enumerate}
 
 \subsection{The inductive step}
 
 So now let us suppose that for $0\ne r\ne s\ne 0$ it is true that $\omega(r,s)=\omega(s,rs)=\omega(rs,r)$ and try to establish the four inductive steps
 
 \begin{enumerate}
 	\item $p=2r,q=2s$\\
 	We will attempt to show that $\omega(2r,2s)=\omega(2s,2rs)=\omega(2rs,2r).$\\
 	$\omega(2r,2s)=-\omega(r,s)$ for $\omega_{0}$ through $\omega_{3}$ and
 	$\omega_{24}$ through $\omega_{27}.$ 
 	
 	$\omega(2r,2s)=\omega(r,s)$ for $\omega_{4}$ through $\omega_{7}$ and
 	$\omega_{28}$ through $\omega_{31}.$ \\
 	
 	$\omega(2s,2rs)=-\omega(s,rs)$ for $\omega_{0}$ through $\omega_{3}$ and
 	$\omega_{24}$ through $\omega_{27}.$ 
 	$\omega(2s,2rs)=\omega(s,rs)$ for $\omega_{4}$ through $\omega_{7}$ and
 	$\omega_{28}$ through $\omega_{31}.$ \\
 
 	$\omega(2rs,2r)=\omega(rs,r)$ for $\omega_{0}$ through $\omega_{3}$ and
 	$\omega_{24}$ through $\omega_{27}.$ 
 	$\omega(2rs,2r)=\omega(rs,r)$ for $\omega_{4}$ through $\omega_{7}$ and
 	$\omega_{28}$ through $\omega_{31}.$ \\
 	
 	Thus $\omega(2r,2s)=\omega(2s,2rs)=\omega(2rs,2r)$ for $\omega_{0}$ through $\omega_{7}$
 	and $\omega_{24}$ through $\omega_{31}.$
 	
 	\item $p=2r,q=2s+1$\\
 	We will attempt to show that $\omega(2r,2s+1)=\omega(2s+1,2rs+1)=\omega(2rs+1,2r).$\\
 	$\omega(2r,2s+1)=\omega(r,s)$ for $\omega_{0},\omega_{1},\omega_{4},\omega_{5},\omega_{24},\omega_{25},\omega_{28},\omega_{29}.$\\
 	$\omega(2r,2s+1)=-\omega(r,s)$ for $\omega_{2},\omega_{3},\omega_{6},\omega_{7},\omega_{26},\omega_{27},\omega_{30},\omega_{31}.$\\
  	
 	$\omega(2s+1,2rs+1)=\omega(s,rs)$ for $\omega_{0},\omega_{2},\omega_{4},\omega_{6},\omega_{24},\omega_{26},\omega_{28},\omega_{30}.$\\
 	$\omega(2r,2s+1)=-\omega(s,rs)$ for $\omega_{1},\omega_{3},\omega_{5},\omega_{7},\omega_{25},\omega_{27},\omega_{29},\omega_{31}.$\\
 	
 	So this property fails for $\omega_{1},\omega_{2},\omega_{5},\omega_{6},\omega_{25},\omega_{26},\omega_{29},\omega_{30}.$
 	
 	$\omega(2rs+1,2r)=\omega(r,s)$ for $\omega_{0},\omega_{4},\omega_{24},\omega_{28}.$\\
 	$\omega(2rs+1,2r)=-\omega(r,s)$ for $\omega_{3},\omega_{7},\omega_{27},\omega_{31}.$
 	
 	Thus $\omega(2r,2s+1)=\omega(2s+1,2rs+1)=\omega(2rs+1,2r)$ for $\omega_{0},\omega_{3},\omega_{4},\omega_{7},\omega_{24},\omega_{27},\omega_{28},\omega_{31}.$
 	
 	\item $p=2r+1,q=2s$\\
 	The result $\omega(2r+1,2s)=\omega(2s+1,2rs+1)=\omega(2rs+1,2r)$ follows immediately from  $\omega(2r,2s+1)=\omega(2s+1,2rs+1)=\omega(2rs+1,2r)$ and the fact that $\omega(p,q)=-\omega(q,p)$ for $0\ne p\ne q\ne 0.$
 	
 	\item $p=2r+1,q=2r+1$\\
 	We will attempt to show that $\omega(2r+1,2s+1)=\omega(2s+1,2rs)=\omega(2rs,2r+1)$ for $\omega_{0},\omega_{3},\omega_{4},\omega_{7},\omega_{24},\omega_{27},\omega_{28},\omega_{31}.$
 	
 	 $\omega(2r+1,2s+1)=\omega(r,s)$ for $\omega_{0},\omega_{4},\omega_{24},\omega_{28}.$\\
 	 $\omega(2r+1,2s+1)=-\omega(r,s)$ for $\omega_{3},\omega_{7},\omega_{27},\omega_{31}.$\\
 	 
 	 $\omega(2s+1,2rs)=\omega(s,rs)$ for $\omega_{0},\omega_{4},\omega_{24},\omega_{28}.$\\
 	$\omega(2s+1,2rs)=\omega(s,rs)$ for $\omega_{3},\omega_{7},\omega_{27},\omega_{31}.$\\
 	
  	 $\omega(2rs,2r+1)=\omega(rs,r)$ for $\omega_{0},\omega_{4},\omega_{24},\omega_{28}.$\\
 	$\omega(2rs,2r+1)=-\omega(rs,r)$ for $\omega_{3},\omega_{7},\omega_{27},\omega_{31}.$\\
 		
 \end{enumerate}
 
 Thus the quaternion properties are satisfied only by the eight twists
$\omega_{0},\omega_{3},\omega_{4},\omega_{7},\omega_{24},\omega_{27},\omega_{28},\omega_{31}.$ 
 Their twist constants are given in Table \vref{Tfifth}. We notice that what distinguishes these eight from the remaining 24 is that
$\alpha^\prime=-\alpha$ and that $\gamma^\prime=\gamma.$ Only in such cases will the product satisfy the quaternion properties. 

\newpage

Table \vref{Ttwist8} shows all eight twist tables
for the octonions. To obtain the traditional quaternion sub-table, one must let $\alpha=1.$
 
\begin{table}[ht]
\begin{tabular}{lrrrrr}
& $\alpha^\prime$ & $\alpha$ & $\beta$ & $\gamma$ & $\gamma^\prime$\\\toprule
$P_{0}:   (a,b)(c,d)=(ca-b^*d,da^*+bc)$&$-$ &$+$ &$-$ &$+$ & $+$ \\
$P_{3}:   (a,b)(c,d)=(ca-db^*,a^*d+cb)$&$-$ &$+$ &$-$ &$-$&$-$ \\
$P_{4}:   (a,b)(c,d)=(ac-b^*d,da^*+bc)$&$-$ &$+$ &$+$ &$+$ &$+$ \\
$P_{7}:   (a,b)(c,d)=(ac-db^*,a^*d+cb)$&$-$ &$+$ &$+$ &$-$ &$-$ \\
$P_{24}:  (a,b)(c,d)=(ca-bd^*,ad+c^*b)$&$+$ &$-$ &$-$&$+$ &$+$ \\
$P_{27}:  (a,b)(c,d)=(ca-d^*b,da+bc^*)$&$+$ &$-$ &$-$&$-$ &$-$ \\
$P_{28}:  (a,b)(c,d)=(ac-bd^*,ad+c^*b)$&$+$ &$-$ &$+$ &$+$ &$+$ \\
$P_{31}:  (a,b)(c,d)=(ac-d^*b,da+bc^*)$&$+$ &$-$ &$+$ &$-$ &$-$ \\\bottomrule
\end{tabular}\\[12pt]
\caption{Twist Constants for the Cayley-Dickson Products}\label{Tfifth}
\end{table}

\begin{table}[ht]
~\\
\Tree [.$\mathbf{C}$ [. $\mathbf{C}$ $\mathbf{T}$ ] [. $\mathbf{L}$ $\mathbf{-D}$ ] ]
\Tree [.$\mathbf{L}$ [. $\mathbf{L}$ $-\alpha\mathbf{N}$ ] [. $\mathbf{L}$ $\alpha\mathbf{N}$ ] ]\\[12pt]
\Tree [.$\mathbf{T}$ [. $\mathbf{T}$ $\mathbf{T}$ ] [. $\alpha\mathbf{N}$ $-\alpha\mathbf{N}$ ] ]
\Tree [.$\mathbf{-D}$ [. $\mathbf{-D}$ $\alpha\mathbf{N}$ ] [ $-\alpha\mathbf{N}$ $\mathbf{-D}$ ] ]\\[12pt]
\Tree [.$\mathbf{N}$ [. $\beta\mathbf{N}$ $\gamma\mathbf{N}$ ] [. $\gamma\mathbf{N}$ $\gamma\mathbf{N}$ ] ]
\Tree [.$-\mathbf{N}$ [. $-\beta\mathbf{N}$ $-\gamma\mathbf{N}$ ] [. $-\gamma\mathbf{N}$ $-\gamma\mathbf{N}$ ] ] 
\vspace{3mm}
\caption{Generalized Cayley-Dickson (actual) Twist Tree}
\label{fig:cydtree}
\end{table}

\begin{table}[ht]
\begin{tabular}{l|rr|rr|rr|rr}
  & 0 &  1        & 2              & 3                & 4               & 5                 & 6               & 7               \\\midrule
0 & 1 &  1        & 1              & 1                & 1               & 1                 & 1               & 1               \\
1 & 1 &  $-1$     & $\alpha$      & $-\alpha$         & $\alpha$       & $-\alpha$          & $\alpha$       & $-\alpha$         \\\midrule
2 & 1 &  $-\alpha$ & $-1$           & $\alpha$        & $\alpha\beta$ & $\alpha\gamma$   & $-\alpha\beta$  & $-\alpha\gamma$   \\
3 & 1 &  $\alpha$  & $-\alpha$       & $-1$             & $\alpha\gamma$ & $\alpha\gamma$ & $-\alpha\gamma$  & $-\alpha\gamma$ \\\midrule
4 & 1 &  $-\alpha$ & $-\alpha\beta$ & $-\alpha\gamma$   & $-1$            & $\alpha$         & $\alpha\beta$ & $\alpha\gamma$   \\
5 & 1 &  $\alpha$  & $-\alpha\gamma$ & $-\alpha\gamma$ & $-\alpha$        & $-1$              & $\alpha\gamma$ & $\alpha\gamma$ \\\midrule
6 & 1 &  $-\alpha$ & $\alpha\beta$  & $\alpha\gamma$    & $-\alpha\beta$  & $-\alpha\gamma$    & $-1$            & $\alpha$         \\
7 & 1 &  $\alpha$  & $\alpha\gamma$  & $\alpha\gamma$  & $-\alpha\gamma$  & $-\alpha\gamma$  & $-\alpha$        & $-1$
\end{tabular}\\[12pt]
  \caption{Octonion Twist Table for the eight Cayley-Dickson products}
  \label{Ttwist8}
\end{table}

\begin{table}[ht]
\begin{tabular}{l|rr|rr|rr|rr}
  & 0 &  1        & 2              & 3                & 4               & 5                 & 6               & 7               \\\midrule
0 & 1 &  1        & 1              & 1                & 1               & 1                 & 1               & 1               \\
1 & 1 &  $-1$     & $1$      & $-1$         & $1$       & $-1$          & $1$       & $-1$         \\\midrule
2 & 1 &  $-1$ & $-1$           & $1$        & $\beta$ & $\gamma$   & $-\beta$  & $-\gamma$   \\
3 & 1 &  $1$  & $-1$       & $-1$             & $\gamma$ & $\gamma$ & $-\gamma$  & $-\gamma$ \\\midrule
4 & 1 &  $-1$ & $-\beta$ & $-\gamma$   & $-1$            & $1$         & $\beta$ & $\gamma$   \\
5 & 1 &  $1$  & $-\gamma$ & $-\gamma$ & $-1$        & $-1$              & $\gamma$ & $\gamma$ \\\midrule
6 & 1 &  $-1$ & $\beta$  & $\gamma$    & $-\beta$  & $-\gamma$    & $-1$            & $1$         \\
7 & 1 &  $1$  & $\gamma$  & $\gamma$  & $-\gamma$  & $-\gamma$  & $-1$        & $-1$
\end{tabular}\\[12pt]
  \caption{Octonion Twist Table when $\alpha=1$}
  \label{Ttwist4}
\end{table}

\begin{table}

\begin{tabular}{ccccc}
$\mathbf{C}$  & $\mathbf{L}$ & $\mathbf{T}$   & $\mathbf{D}$   & $\mathbf{N}$ \\\toprule
$r=s=0$ & $r>s=0$   & $s>r=0$    & $r=s\ne0$  & $ 0\ne r\ne s\ne 0$ \\\midrule
$\left[ \begin{array}{rr}  1 &  1\\  1 & -1 \end{array} \right]$
                                & $\left[ \begin{array}{rr}  1 & -1\\  1 &  1 \end{array} \right]$
                               & $\left[ \begin{array}{rr}  1 &  1\\  1 & -1 \end{array} \right]$
                               & $\left[ \begin{array}{rr}  1 & -1\\  1 &  1 \end{array} \right]$
                               & $\left[ \begin{array}{rr} \beta &  \gamma\\  \gamma &  \gamma \end{array} \right]$\\\midrule
\end{tabular}\\[12pt]
\caption{Twist Table Components for $\omega_{0},\omega_{3},\omega_{4},\omega_{7}$}
\end{table}

\section{The product $P_{0}$ and Octonion Index Cycling}

In Cayley-Dickson spaces, if $i_pi_q=i_r$ then that fact is abbreviated $(p,q,r).$ This is called a structure constant.

If $0\ne p\ne q\ne 0$ then for all 32 products one and only one of $(p,q,pq)$ or $(q,p,pq)$ is true. For those eight which satisfy the
quaternion property, $(p,q,pq)$ implies $(q,pq,p)$ which implies $(pq,p,q.)$

If $p$ is an integer and $1\le p\le 7$ and $\pi_0$ is the permutation (1376524), then define the \emph{$\pi_0$ successor} $p^\prime$ of $p$ as the
integer following $p$ in the permutation. Then for the  product $P_{0}$ alone it is true that $(p,q,r)$ implies $(p^\prime,q^\prime,r^\prime).$ 

\begin{eqnarray*} (1,2,3)&\to(3,4,7)\\&\to(7,1,6)\\&\to(6,3,5)\\&\to(5,7,2)\\&\to(2,6,4)\\&\to(4,5,1)\\&\to(1,2,3)
\end{eqnarray*}
This, together with the second quaternion property is sufficient to recover the $P_{0}$ multiplication table for the octonions.

\begin{table}[ht]
\begin{center}
\begin{tabular}{|r|r|r|r|r|r|r|r|}
\hline
$1$ &
$i_1$ &
$i_2$ &
$i_3$ &
$i_4$ &
$i_5$ &
$i_6$ &
$i_7$\\\hline
$i_1$ &
$-1$ &
$i_3$ &
$-i_2$ &
$i_5$ &
$-i_4$ &
$i_7$ &
$-i_6$\\\hline
$i_2$ &
$-i_3$ &
$-1$ &
$i_1$ &
$-i_6$ &
$i_7$ &
$i_4$ &
$-i_5$\\\hline
$i_3$ &
$i_2$ &
$-i_1$ &
$-1$ &
$i_7$ &
$i_6$ &
$-i_5$ &
$i_4$\\\hline
$i_4$ &
$-i_5$ &
$i_6$ &
$-i_7$ &
$-1$ &
$i_1$ &
$-i_2$ &
$i_3$\\\hline
$i_5$ &
$i_4$ &
$-i_7$ &
$-i_6$ &
$-i_1$ &
$-1$ &
$i_3$ &
$i_2$\\\hline
$i_6$ &
$-i_7$ &
$-i_4$ &
$i_5$ &
$i_2$ &
$-i_3$ &
$-1$ &
$i_1$\\\hline
$i_7$ &
$i_6$ &
$i_5$ &
$i_4$ &
$-i_3$ &
$-i_2$ &
$-i_1$ &
$-1$\\\hline
\end{tabular}\\[12pt]

\caption{Basis Vector Multiplication Table for $P_0$}
\end{center}
\end{table}

\section{Recursive definition of structure constants for the Cayley-Dickson twists}

For $\omega_{0},\omega_{3},\omega_{4}$ and $\omega_{7},$ 
 
\begin{equation}(1,2n,2n+1) \text{\ for all\ } n>0\end{equation}

whereas for  $\omega_{24},\omega_{27},\omega_{28}$ and $\omega_{31},$ 

\begin{equation}(2n,1,2n+1) \text{\ for all\ } n>0\end{equation}

And for all eight of the Cayley-Dickson twists, the second quaternion property holds.

For $0\ne p\ne q\ne 0$

\begin{equation}(p,q,r)\longrightarrow (q,r,p)\longrightarrow (r,p,q)\end{equation}

For $\omega_{0}$ and $\omega_{24}$ if $0\ne p\ne q\ne 0$ then 
\begin{align}
      (p,q,r) &\longrightarrow (2q,2p,2r)\notag\\
                    &\longrightarrow (2p,2q+1,2r+1)\notag\\
                    &\longrightarrow (2p+1,2q,2r+1)\notag\\
                    &\longrightarrow (2p+1,2q+1,2r)
\end{align}

For $\omega_{3}$ and $\omega_{27}$ if $0\ne p\ne q\ne 0$ then 

\begin{align}
      (p,q,r) &\longrightarrow (2q,2p,2r)\notag\\
                    &\longrightarrow (2q,2p+1,2r+1)\notag\\
                    &\longrightarrow (2q+1,2p,2r+1)\notag\\
                    &\longrightarrow (2q+1,2p+1,2r)
\end{align}

For $\omega_{4}$ and $\omega_{25}$ if $0\ne p\ne q\ne 0$ then 
\begin{align}
      (p,q,r) &\longrightarrow (2p,2q,2r)\notag\\
                    &\longrightarrow (2p,2q+1,2r+1)\notag\\
                    &\longrightarrow (2+1,2r,2r+1)\notag\\
                    &\longrightarrow (2p+1,2q+1,2r)
\end{align}

For $\omega_{7}$ and $\omega_{31}$ if $0\ne p\ne q\ne 0$ then 
\begin{align}
      (p,q,r) &\longrightarrow (2p,2q,2r)\notag\\
                    &\longrightarrow (2q,2p+1,2r+1)\notag\\
                    &\longrightarrow (2q+1,2p,2r+1)\notag\\
                    &\longrightarrow (2q+1,2p+1,2r)
\end{align}

\section{Octonion Cycles}

If $0\ne p\ne q\ne 0$, and if $i_pi_q=+i_r$ we say that the \emph{sense} of $(p,q,r)$ (and therefore of $(q,r,p)$ and $r,p,q)$) is \emph{positive} or
that $\sigma(p,q,r)=1.$ Otherwise we say that the the sense is \emph{negative} and $\sigma(p,q,r)=-1.$

So it follows that for the twist constants $\alpha,\beta$ and $\gamma$

\begin{equation}
  \alpha = \sigma(1,2,3)
\end{equation}

\begin{equation}
  \alpha\beta = \sigma(2,4,6)
\end{equation}

\begin{equation}
  \alpha\gamma = \sigma(2,5,7)
\end{equation}

If we define $\pi$ as the permutation $(1)(357)(246)$ and for $1\le p\le 7$ we
let $p^{\prime\prime}$ denote the $\pi$ successor of $p$, then

\begin{equation}
  \sigma(p,q,r)=\sigma(p^{\prime\prime},q^{\prime\prime},r^{\prime\prime})
\end{equation}

Thus, for any of the eight Cayley-Dickson doubling variations, knowing the sense
of $(1,2,3),$ $(2,4,6)$ and $(2,5,7)$ is sufficient to recover the sense of all 42
triples $(p,q,r).$

\begin{figure}
 \psset{xunit=1.0cm,yunit=1.0cm,algebraic=true,dotstyle=o,dotsize=3pt 0,
linewidth=1.5pt,arrowsize=3pt 2,arrowinset=0.25}
\begin{pspicture*}(-4.3,-3.81)(30.68,6.3)
\pspolygon[linestyle=none,fillstyle=solid,fillcolor=white](-3.36,-2.8)(5.96,-2.92)(1.41,5.21)
\psline[linecolor=red](-3.36,-2.8)(5.96,-2.92)
\psline[linecolor=red](5.96,-2.92)(1.41,5.21)
\psline[linecolor=red](1.41,5.21)(-3.36,-2.8)
\psline[linecolor=green](5.96,-2.92)(-0.98,1.2)
\psline[linecolor=green](1.41,5.21)(1.3,-2.86)
\psline[linecolor=green](-3.36,-2.8)(3.69,1.14)
\pscircle[linecolor=blue](1.34,-0.17){2.69}
\psdots[dotstyle=*,linecolor=black](-3.36,-2.8)
\rput[bl](-3.92,-3.09){\black{$\mathbf{i_7}$}}
\psdots[dotstyle=*,linecolor=black](5.96,-2.92)
\rput[bl](6.18,-3.23){\black{$\mathbf{i_5}$}}
\psdots[dotstyle=*,linecolor=black](1.41,5.21)
\rput[bl](1.34,5.56){\black{$\mathbf{i_3}$}}
\psdots[dotstyle=*,linecolor=black](3.69,1.14)
\rput[bl](3.8,1.32){\black{$\mathbf{i_6}$}}
\psdots[dotstyle=*,linecolor=black](1.3,-2.86)
\rput[bl](1.23,-3.55){\black{$\mathbf{i_2}$}}
\psdots[dotstyle=*,linecolor=black](-0.98,1.2)
\rput[bl](-1.45,1.4){\black{$\mathbf{i_4}$}}
\psdots[dotstyle=*,linecolor=black](1.34,-0.17)
\rput[bl](1.48,0.25){\black{$\mathbf{i_1}$}}
\end{pspicture*}
\caption{Fano Plane}
\label{fano}
\end{figure}
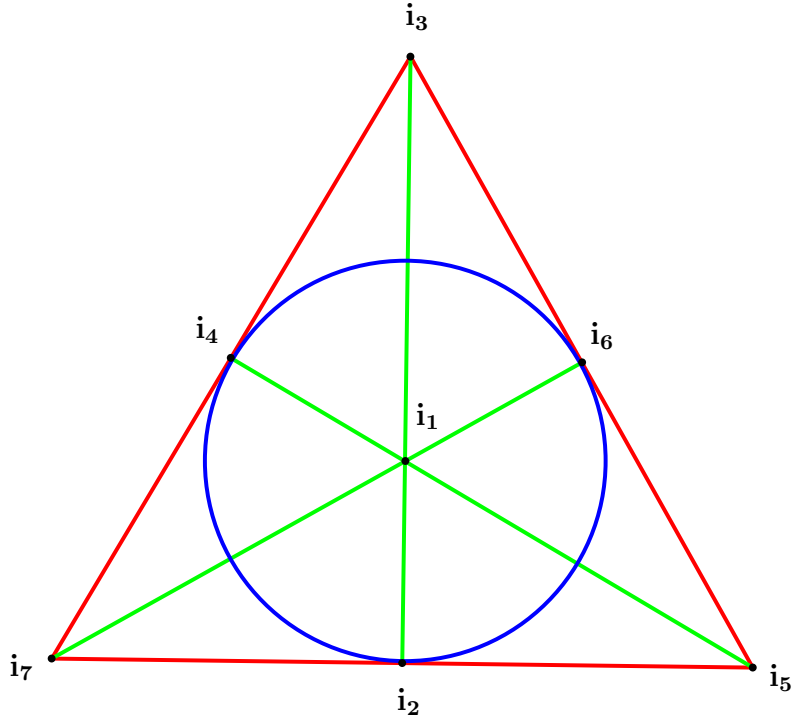

One way to recover any of the eight variations on the sets of 42 triples is by
using the Fano plane (Figure \vref{fano}).

Each of the three sides of the triangle represents a triple $(p,q,r).$ So does
each altitude from a vertex to the midpoint of the opposite side and so does the
circle through the midpoints. For each of the eight Cayley-Dickson products, the
sense of the three sides is the same--either clockwise ($\leftarrow$) around
the triangle, or counter-clockwise ($\rightarrow$). If clockwise, then the three
sides of the triangle represent $(5,2,7),$ $(7,4,3)$ and $(3,6,5).$ If
counter-clockwise then the three sides represent $(7,2,5),$ $(5,6,3)$ and
$(3,4,7).$ The circle through the midpoints of the sides represents either
$(2,4,6)$ in the clockwise sense ($\circlearrowright$) or $(2,6,4)$ in
the counter-clockwise sense ($\circlearrowleft$). The three altitudes may all
be in an `up' sense from center to vertex ($\uparrow$) or
may all be in a `down' sense from center to base ($\downarrow$).
So the altitudes must either be $(2,1,3),$ $(6,1,7),$ and $(4,1,5)$ or $(3,1,2),$ $(7,1,6),$ and $(5,1,4)$.
All altitudes must have the same sense. See Table
\vref{Tfano} for a breakdown of all the modes of the eight Cayley-Dickson
doubling products. A thing to notice about the Fano plane is that a $120^\circ$
counter-clockwise rotation of the diagram represents the permutation
$\pi=(1)(357)(246)$ mentioned above.

Since the sides may have two senses and the circle may have two senses and the
altitudes may have two senses, all $2^3=8$ versions of the Cayley-Dickson
products may be accommodated in the one diagram.

\begin{table}
\begin{tabular}{lrrr}
&&&\\\toprule
$P_{0}:   (a,b)(c,d)=(ca-b^*d,da^*+bc)$&$\downarrow$ &$\circlearrowleft$
&$\rightarrow$\\
$P_{3}:   (a,b)(c,d)=(ca-db^*,a^*d+cb)$&$\downarrow $ &$\circlearrowleft$
&$\leftarrow$\\
$P_{4}:   (a,b)(c,d)=(ac-b^*d,da^*+bc)$&$\downarrow $ &$\circlearrowright$
&$\rightarrow$\\
$P_{7}:   (a,b)(c,d)=(ac-db^*,a^*d+cb)$&$\downarrow $ &$\circlearrowright$
&$\leftarrow$\\
$P_{24}:  (a,b)(c,d)=(ca-bd^*,ad+c^*b)$&$\uparrow$ &$\circlearrowright$
&$\leftarrow$\\
$P_{27}:  (a,b)(c,d)=(ca-d^*b,da+bc^*)$&$\uparrow$ &$\circlearrowright$
&$\rightarrow$\\
$P_{28}:  (a,b)(c,d)=(ac-bd^*,ad+c^*b)$&$\uparrow $ &$\circlearrowleft$
&$\leftarrow$\\
$P_{31}:  (a,b)(c,d)=(ac-d^*b,da+bc^*)$&$\uparrow $ &$\circlearrowleft$
&$\rightarrow$\\\bottomrule
\end{tabular}\\[12pt]
\caption{Fano Plane Modes}\label{Tfano}
\end{table}

\end{document}